\documentclass[11pt,reqno]{amsart}
\usepackage[T1]{fontenc}
\usepackage{lmodern}
\usepackage[utf8]{inputenc}
\usepackage{amsmath,amssymb,mathtools,bm}
\usepackage[margin=1.05in,top=1in,bottom=1in]{geometry}
\usepackage{microtype}
\usepackage[hidelinks]{hyperref}
\hypersetup{pdftitle={The BCC Lattice Maximizes the First Positive Neumann Eigenvalue among 3D Lattice Voronoi Cells},pdfauthor={Kaixin Deng and Senping Luo},pdfsubject={Neumann eigenvalues and lattice Voronoi cells}}
\numberwithin{equation}{section}
\newtheorem{theorem}{Theorem}[section]
\newtheorem{proposition}[theorem]{Proposition}
\newtheorem{lemma}[theorem]{Lemma}

\newtheorem{numlemma}{Numerical Lemma}
\theoremstyle{remark}
\newtheorem{remark}[theorem]{Remark}
\newcommand{\R}{\mathbb R}
\newcommand{\Z}{\mathbb Z}
\newcommand{\Nzero}{\mathbb N_0}
\newcommand{\BCC}{\mathrm{BCC}}
\newcommand{\dd}{\,\mathrm d}
\DeclareMathOperator{\covol}{covol}
\DeclareMathOperator{\diag}{diag}
\DeclareMathOperator{\off}{off}
\DeclareMathOperator{\tr}{tr}
\DeclareMathOperator{\adj}{adj}
\DeclareMathOperator{\sgn}{sgn}

\DeclareMathOperator{\argmax}{arg\,max}
\newcommand{\norm}[1]{\lVert#1\rVert}
\newcommand{\abs}[1]{\lvert#1\rvert}
\newcommand{\1}{\mathbf 1}
\allowdisplaybreaks[2]
\title[BCC and the first Neumann eigenvalue]{The BCC Lattice Maximizes the First Positive Neumann Eigenvalue among 3D Lattice Voronoi Cells}
\author{Kaixin Deng}
\author{Senping Luo}

\address[K.~Deng]{School of Mathematics and statistics, Jiangxi Normal University, Nanchang, 330022, China}

\email[K.~Deng]{Dengkaixin1999@163.com}

\address[S.~Luo]{School of Mathematics and statistics, Jiangxi Normal University, Nanchang, 330022, China}

\email[S.~Luo]{luosp1989@163.com}

\subjclass[2020]{Primary 35P15; Secondary 11H06, 52C07}
\keywords{Neumann eigenvalue, lattice Voronoi cell, body-centered cubic lattice, spectral optimization, Selling reduction}
\begin{document}
\begin{abstract}
It is shown that, among three-dimensional lattices of fixed covolume, the body-centered cubic lattice maximizes the first positive Neumann eigenvalue of its Voronoi cell uniquely, up to orthogonal transformations.
\end{abstract}
\maketitle

\section{Introduction}\label{sec:intro}
We prove that BCC uniquely maximizes the first positive Neumann eigenvalue at fixed covolume by Selling reduction, local eigenvalue splitting and global Rayleigh bounds.

Let $\Lambda\subset\R^3$ be a full-rank lattice, and let
\begin{equation}\label{eq:voronoi}
 V_\Lambda=\{x\in\R^3:\abs{x}\le\abs{x-v}\ \text{for all }v\in\Lambda\}
\end{equation}
be its Voronoi cell. For a bounded connected Lipschitz domain $\Omega$, write $\mu_1(\Omega)$ for the first \emph{positive} Neumann eigenvalue:
\begin{equation}\label{eq:rayleigh}
 \mu_1(\Omega)=
 \inf_{\substack{u\in H^1(\Omega)\setminus\{0\}\\ \int_\Omega u=0}}
 \frac{\int_\Omega\abs{\nabla u}^2}{\int_\Omega u^2}.
\end{equation}
Since
\begin{equation}\label{eq:scaling}
 \abs{V_\Lambda}=\covol\Lambda,
 \qquad \mu_1(a\Omega)=a^{-2}\mu_1(\Omega),
\end{equation}
we normalize $\covol\Lambda=1$.

The Szeg\H{o}--Weinberger inequality \cite{Szego,Weinberger} gives
\[
 |\Omega|^{2/3}\mu_1(\Omega)\le |B|^{2/3}\mu_1(B)
\]
where $B\subset\R^3$ is a ball; equality holds precisely when $\Omega$ is a ball.
Brasco and Pratelli \cite{BrascoPratelli} established sharp quantitative stability. We determine the maximizer among lattice Voronoi cells. Set
the body-centered cubic lattice as follows:
\begin{equation}\label{eq:bcc}
 \BCC=2^{1/3}\left[\Z^3\cup
 \left(\Z^3+\tfrac12(1,1,1)\right)\right],\qquad \covol\BCC=1.
\end{equation}

\begin{theorem}\label{thm:main}
For every full-rank lattice $\Lambda\subset\R^3$ of covolume one,
\begin{equation}\label{eq:main}
 \mu_1(V_\Lambda)\le\mu_1(V_\BCC),
\end{equation}
with equality if and only if $\Lambda=R\BCC$ for some $R\in O(3)$. The maximizing cell is the regular truncated octahedron.
\end{theorem}

\begin{samepage}
\begin{remark}[Covolume and density]\label{rem:scaling}
For covolume $v>0$ and density $\varrho=v^{-1}$, Theorem~\ref{thm:main} gives
\begin{equation}\label{eq:fixed-density}
 \max_{\covol\Lambda=v}\mu_1(V_\Lambda)
 =v^{-2/3}\mu_1(V_\BCC)
 =\varrho^{2/3}\mu_1(V_\BCC).
\end{equation}
Without prescribing the scale,
\[
 \argmax_{\Lambda\subset\R^3}
 \bigl[(\covol\Lambda)^{2/3}\mu_1(V_\Lambda)\bigr]
 =\{aR\BCC:a>0,\ R\in O(3)\}.
\]
\end{remark}
\end{samepage}

Numerically, to ten decimal places,
\begin{equation}\label{eq:bcc-value}
 \mu_1(V_\BCC)=10.8374401273.
\end{equation}

For the insulated heat equation,
\[
 \left.
 \begin{gathered}
 \partial_tu=\Delta u,\qquad \partial_\nu u=0,\\
 \int_\Omega u(0)=0
 \end{gathered}\right\}
 \quad\Longrightarrow\quad
 \norm{u(t)}_{L^2(\Omega)}
 \le e^{-\mu_1(\Omega)t}\norm{u(0)}_{L^2(\Omega)}.
\]
The decay rate is sharp; Theorem~\ref{thm:main} gives the fastest decay among lattice cells of fixed volume.

For unit-area polygons, Laugesen and Siudeja
\cite{LaugesenSiudeja2009} and Endo and Osting \cite{EndoOsting} proved,
respectively,
\[
 \max_{\substack{\Omega\text{ a triangle}\\|\Omega|=1}}\mu_1(\Omega)
 =\frac{4\pi^2}{3\sqrt3},\qquad
 \max_{\substack{\Omega\text{ a convex quadrilateral}\\|\Omega|=1}}\mu_1(\Omega)
 =\pi^2.
\]
The maximizers are the equilateral triangle and the square.
For higher-dimensional eigenvalue sums, see \cite{LaugesenSiudeja2011}.

BCC also minimizes the Voronoi moments:
\[
 \int_{V_\Lambda}|x|^p\dd x\ge\int_{V_\BCC}|x|^p\dd x,
 \qquad \covol\Lambda=1,\quad 1\le p<\infty.
\]
Barnes and Sloane \cite{BarnesSloane} proved the case $p=2$; the second author
\cite{LuoQuantizer} proved every finite $p\ge1$, resolving the lattice
quantizer conjecture of Graf and Luschgy \cite{GrafLuschgy}.
For quantization and centroidal Voronoi tessellations, see
\cite{DuFaberGunzburger,GrafLuschgy}; for surface-area optimization in
Selling coordinates, see \cite{CesaroniNovaga}.

Section~\ref{sec:reduction} gives the conorm representation and geometric bounds; Section~\ref{sec:spectrum} identifies the BCC eigenspace and encloses its eigenvalue. Sections~\ref{sec:local} and \ref{sec:global} establish the local and global comparisons. The numerical computations are put in the Appendix.

\section{Selling reduction and geometric bounds}\label{sec:reduction}
Six nonnegative Selling conorms parametrize lattice shapes. Dual sine trials and projected moment comparisons give the short-vector bounds.

Put $I=\{0,1,2,3\}$ and list the six edges as
\[
 E=\{01,02,03,12,13,23\},\qquad ij=\{i,j\}\quad(0\le i<j\le3).
\]
We write $p=(p_{01},p_{02},p_{03},p_{12},p_{13},p_{23})$ and set
$p_{ji}=p_{ij}$ for $i\ne j$. A superbase $(b_i)_{i=0}^3$ satisfies $\sum_{i=0}^3b_i=0$, with any three vectors forming a lattice basis. Every three-dimensional lattice has an obtuse superbase \cite{ConwaySloane}. Indeed, if $b_i\cdot b_j>0$, the Selling move
\[
 (b_i,b_j,b_k,b_l)\longmapsto
 (-b_i,b_j,b_k+b_i,b_l+b_i)
\]
preserves the superbase property and satisfies
\[
 \sum_{r=0}^3\abs{b'_r}^2-\sum_{r=0}^3\abs{b_r}^2
 =2b_i\cdot(b_k+b_l)+2\abs{b_i}^2
 =-2b_i\cdot b_j<0.
\]
Only finitely many superbases have bounded squared length, so the reduction terminates.

Normalize the conorms by
\begin{equation}\label{eq:conorms}
 p_{ij}=\frac{-b_i\cdot b_j}{\sum_{0\le a<b\le3}(-b_a\cdot b_b)},
 \qquad p\in\Delta:=\left\{p\ge0:p_{01}+p_{02}+p_{03}+p_{12}+p_{13}+p_{23}=1\right\}.
\end{equation}
Put $e_0=0$, let $e_1,e_2,e_3$ be the coordinate vectors, and define
\begin{equation}\label{eq:gram}
 \begin{gathered}
 d_{ij}=e_i-e_j,\qquad
 G(p)=\sum_{0\le i<j\le3}p_{ij}d_{ij}d_{ij}^T,\\
 D(p)=\det G(p),\qquad
 C_{ij}(p)=d_{ij}^T\adj(G(p))d_{ij}\quad(0\le i<j\le3).
 \end{gathered}
\end{equation}
For $e=ij\in E$, the edge notation means
$p_e=p_{ij}$, $d_e=d_{ij}$, and $C_e=C_{ij}$. In particular,
\[
 G(p)=\begin{pmatrix}
 p_{01}+p_{12}+p_{13}&-p_{12}&-p_{13}\\
 -p_{12}&p_{02}+p_{12}+p_{23}&-p_{23}\\
 -p_{13}&-p_{23}&p_{03}+p_{13}+p_{23}
 \end{pmatrix}.
\]
The identity
\[
 \left|\sum_{i=1}^3 n_i b_i\right|^2
 =\sum_{0\le i<j\le3}(-b_i\cdot b_j)(n_i-n_j)^2,
 \qquad n_0=0,
\]
gives
\[
 \operatorname{Gram}(b_1,b_2,b_3)
 =\left(\sum_{0\le a<b\le3}(-b_a\cdot b_b)\right)G(p),
 \qquad \det(D^{-1/6}G^{1/2})=1.
\]
Thus
\begin{equation}\label{eq:lattice-p}
 \Lambda(p)=D(p)^{-1/6}G(p)^{1/2}\Z^3,
 \qquad p\in\Delta,\quad D(p)>0,
\end{equation}
parametrizes unit-covolume lattices up to orthogonal equivalence. For $n_0=0$,
\[
 n^TGn=0
 \ \Longleftrightarrow\ n_i=n_j\quad\text{whenever }p_{ij}>0,
\]
and therefore
\[
 D(p)>0
 \ \Longleftrightarrow\ (I,\{ij\in E:p_{ij}>0\})\text{ is connected}.
\]

For $\varnothing\ne S\subsetneq I$, set
\[
 W_S(p)=\sum_{i\in S}\sum_{j\in I\setminus S}p_{ij},\qquad
 Y(p)=\sum_{0\le i<j\le3}[-p_{ij}d_{ij}/2,p_{ij}d_{ij}/2].
\]
With $y_0=-y_1-y_2-y_3$, we have
\begin{align}
 Y(p)&=\left\{y\in\R^3:\sum_{i\in S}y_i\le\tfrac12W_S(p)
                  \quad(\varnothing\ne S\subsetneq I)\right\},\label{eq:zonotope}\\
 V_{\Lambda(p)}&=D^{-1/6}G^{-1/2}Y(p).\label{eq:cell-p}
\end{align}
In $y$ coordinates, the Voronoi inequalities are
\[
 2n\cdot y\le\sum_{0\le i<j\le3}p_{ij}(n_i-n_j)^2,\qquad n\in\Z^3,\quad n_0=0.
\]
For integral $n$,
\[
 \abs{n_i-n_j}\le(n_i-n_j)^2,\qquad
 n_i=\1_S(i)-\1_S(0)\ \Longrightarrow\ n^TGn=W_S.
\]
The segment sum satisfies the Voronoi inequalities; cut indicators recover
the constraints in \eqref{eq:zonotope}. Conversely, for $u_0=0$,
\[
 2u\cdot y
 =2\int_{\min_{0\le i\le3}u_i}^{\max_{0\le i\le3}u_i}\sum_{\substack{0\le i\le3\\u_i>s}}y_i\dd s
 \le\int_{\min_{0\le i\le3}u_i}^{\max_{0\le i\le3}u_i}W_{\{i\in I:u_i>s\}}\dd s
 =\sum_{0\le i<j\le3}p_{ij}\abs{u_i-u_j}.
\]
The last expression is twice the support function of the segment sum,
proving \eqref{eq:zonotope}--\eqref{eq:cell-p}.

By \eqref{eq:cell-p} and \eqref{eq:rayleigh},
\begin{equation}\label{eq:functional}
 F(p):=\mu_1(V_{\Lambda(p)})
 =D^{1/3}\inf_{\substack{f\in H^1(Y(p))\setminus\{0\}\\\int_{Y(p)}f=0}}
 \frac{\int_{Y(p)}\nabla f^TG\nabla f}{\int_{Y(p)}f^2}.
\end{equation}
In particular,
\[
 Y(p)=-Y(p),\quad f(-y)=-f(y)\ \Longrightarrow\ \int_{Y(p)}f=0,
\]
and, for $s>0$,
\[
 G(sp)=sG(p),\quad D(sp)=s^3D(p),\quad
 \Lambda(sp)=\Lambda(p),\quad F(sp)=F(p).
\]

At the equal-conorm point
\begin{equation}\label{eq:equal-point}
 p_*=(1/6,\ldots,1/6),\qquad
 G(p_*)=\frac16\begin{pmatrix}3&-1&-1\\-1&3&-1\\-1&-1&3\end{pmatrix},
 \qquad D(p_*)=\frac2{27}.
\end{equation}
The three permutations of $(1,-1,-1)/2$ form a BCC basis with Gram
matrix $\tfrac32G(p_*)$; hence $\Lambda(p_*)=R\BCC$ for some $R\in O(3)$.
Fix
\begin{equation}\label{eq:thresholds}
 T=\frac{10837}{1000}.
\end{equation}
The proof of Theorem~\ref{thm:main} reduces to
\begin{alignat}{2}
 F(p_*)&>T,&&\qquad\text{Section~\ref{sec:spectrum}},\label{eq:three-a}\\
 0<\norm{p-p_*}_2\le\frac{73}{6000}&\ \Longrightarrow\ F(p)<F(p_*),
 &&\qquad\text{Section~\ref{sec:local}},\label{eq:three-b}\\
 \norm{p-p_*}_2>\frac{73}{6000}&\ \Longrightarrow\ F(p)<T,
 &&\qquad\text{Section~\ref{sec:global}}.\label{eq:three-c}
\end{alignat}

Write $\lambda_1(\Lambda)$ for the shortest nonzero vector length, and let
$\Lambda^*=\{w:w\cdot\Lambda\subset\Z\}$.

\begin{lemma}\label{lem:lengths}
For $p\in\Delta$ with $D(p)>0$,
\begin{equation}\label{eq:lengths}
 \lambda_1(\Lambda(p))^2=D^{-1/3}\min_{\varnothing\ne S\subsetneq I}W_S,
 \qquad
 \lambda_1(\Lambda(p)^*)^2=D^{-2/3}\min_{0\le i<j\le3}C_{ij}.
\end{equation}
\end{lemma}
\begin{proof}
For $n\in\Z^3\setminus\{0\}$, put
\[
 n_0=0,\qquad
 S=\{i\in I:n_i>\min_{0\le j\le3}n_j\},\qquad
 u_i=\1_S(i)-\1_S(0)\quad(1\le i\le3).
\]
Then $\varnothing\ne S\subsetneq I$ and
\[
 n^TGn\ge W_S\ge\min_{\varnothing\ne T\subsetneq I}W_T,
 \qquad u^TGu=W_S.
\]

For the dual lattice, put
\[
 \widetilde m=\left(-\sum_{i=1}^3m_i,m_1,m_2,m_3\right),\qquad
 S=\{i\in I:\widetilde m_i>0\},\qquad
 u_i=\1_S(i)-\1_S(0).
\]
If $m\notin\{0,\pm d_e:e\in E\}$, then
\[
 m\cdot u=\sum_{i\in S}\widetilde m_i\ge2,
 \qquad
 4\le(m\cdot u)^2\le(m^TG^{-1}m)W_S\le m^TG^{-1}m.
\]
On the other hand,
\[
 \sum_{0\le i<j\le3}p_{ij}d_{ij}^TG^{-1}d_{ij}
 =\tr(G^{-1}G)=3,
\]
so
\[
 \min_{m\in\Z^3\setminus\{0\}}m^TG^{-1}m
 =\min_{0\le i<j\le3}d_{ij}^TG^{-1}d_{ij}\le3.
\]
Since
\[
 \Lambda(p)^*=D^{1/6}G^{-1/2}\Z^3,\qquad
 G^{-1}=D^{-1}\adj G,
\]
the scaling factors give \eqref{eq:lengths}.
\end{proof}

Fundamental-cell orthogonality and central symmetry give, for $0\ne w\in\Lambda^*$,
\[
 \int_{V_\Lambda}e^{2\pi i w\cdot x}=0,\qquad
 \int_{V_\Lambda}\sin(\pi w\cdot x)=0,
\]
\[
 \int_{V_\Lambda}\sin^2(\pi w\cdot x)
 =\int_{V_\Lambda}\cos^2(\pi w\cdot x)=\frac{\abs{V_\Lambda}}2.
\]
Consequently,
\[
 \frac{\int_{V_\Lambda}\abs{\nabla\sin(\pi w\cdot x)}^2}
      {\int_{V_\Lambda}\sin^2(\pi w\cdot x)}=\pi^2\abs w^2,
\]
and
\begin{equation}\label{eq:dualbound}
 \mu_1(V_\Lambda)\le\pi^2\lambda_1(\Lambda^*)^2,
 \qquad F(p)\le\pi^2D^{-2/3}\min_{0\le i<j\le3}C_{ij}\le3\pi^2D^{1/3}.
\end{equation}
In particular,
\[
 p\in\Delta,\quad D(p)\longrightarrow0
 \quad\Longrightarrow\quad F(p)\longrightarrow0.
\]

For a primitive $v\in\Lambda$, put
\[
 a=\abs v,\qquad \Gamma=\operatorname{proj}_{v^\perp}\Lambda,
 \qquad \covol\Gamma=a^{-1}.
\]
Let $w(y)$ be the fiber length of $V_\Lambda$ above $y\in v^\perp$.
The bisectors of $\pm v$ and the lattice tiling give
\begin{equation}\label{eq:fibers}
 0\le w(y)\le a,\qquad
 \sum_{\gamma\in\Gamma}w(y+\gamma)=a\quad\text{a.e.}
\end{equation}
For a locally Lipschitz profile $g$ with $g(0)=0$, let
\[
 u_i(y)=g(|y|)\frac{y_i}{|y|}\quad(i=1,2).
\]
Then
\[
 \begin{gathered}
 \sum_{i=1}^2u_i(y)^2=g(|y|)^2,\qquad
 \sum_{i=1}^2|\nabla u_i(y)|^2
 =g'(|y|)^2+\frac{g(|y|)^2}{|y|^2},\\
 \int_{\R^2}u_i(y)w(y)\dd y=0\quad(i=1,2),
 \end{gathered}
\]
since $w(-y)=w(y)$. Suppose $r\mapsto g(r)^2$ is nondecreasing,
$r\mapsto g'(r)^2+g(r)^2/r^2$ is nonincreasing, and
$\int_{V_\Gamma}g(|y|)^2\dd y>0$.
For $y\in V_\Gamma$, $|y+\gamma|\ge|y|$ for every $\gamma\in\Gamma$; hence
\begin{align*}
 &\sum_{\gamma\in\Gamma}
 \left(g'(|y+\gamma|)^2+\frac{g(|y+\gamma|)^2}{|y+\gamma|^2}\right)w(y+\gamma)\\
 &\qquad\le a\left(g'(|y|)^2+\frac{g(|y|)^2}{|y|^2}\right),\\
 &\sum_{\gamma\in\Gamma}g(|y+\gamma|)^2w(y+\gamma)\ge ag(|y|)^2.
\end{align*}
Rayleigh's principle and \eqref{eq:fibers} give
\begin{equation}\label{eq:projected}
 \begin{aligned}
 \mu_1(V_\Lambda)
 &\le\min_{\substack{i=1,2\\\int_{\R^2}u_i^2w>0}}
       \frac{\int_{\R^2}|\nabla u_i|^2w}{\int_{\R^2}u_i^2w}\\
 &\le\frac{\sum_{i=1}^2\int_{\R^2}|\nabla u_i|^2w}
            {\sum_{i=1}^2\int_{\R^2}u_i^2w}\\
 &=\frac{\int_{\R^2}
       \left(g'(|y|)^2+g(|y|)^2/|y|^2\right)w(y)\dd y}
       {\int_{\R^2}g(|y|)^2w(y)\dd y}\\
 &\le\frac{\int_{V_\Gamma}
       \left(g'(|y|)^2+g(|y|)^2/|y|^2\right)\dd y}
       {\int_{V_\Gamma}g(|y|)^2\dd y}.
 \end{aligned}
\end{equation}

\begin{lemma}[Planar disk overlap]\label{lem:disk}
For a planar lattice Voronoi polygon $P$ and the centered regular hexagon $H$ with $\abs H=\abs P$,
\begin{equation}\label{eq:disk}
 \abs{P\cap B_r}\le\abs{H\cap B_r}\qquad(r>0).
\end{equation}
For integrable radial $f$,
\[
 \dd f\ge0\ \Longrightarrow\ \int_Pf\ge\int_Hf,\qquad
 \dd f\le0\ \Longrightarrow\ \int_Pf\le\int_Hf.
\]
\end{lemma}
\begin{proof}
We use the planar moment comparison of \cite{GruberMoment,BoroczkyCsikos}.
A Voronoi facet has relative-interior midpoint $v/2$. If
$w\equiv v\pmod{2\Gamma}$, $w\ne\pm v$, and $\abs w\le\abs v$, then
\[
 \frac{v+w}{2}\in\Gamma\setminus\{0,v\},\qquad
 \left|\frac{v+w}{2}-\frac v2\right|
 =\frac{\abs w}{2}\le\frac{\abs v}{2},
\]
contradicting the relative-interior Voronoi inequalities. In particular,
$w=0$ excludes $v\in2\Gamma$, and
\[
 \#\{\text{facets of }P\}
 \le2\bigl(\#(\Gamma/2\Gamma)-1\bigr)=6.
\]

By scaling, take $r=1$. A facet triangle of area $A$ and half-angle $\theta$ has altitude $\sqrt{A\cot\theta}$ and disk overlap
\[
 H(A,\theta)=\int_0^\theta\min\{1,A\cot\theta\sec^2\varphi\}\dd\varphi.
\]
Equivalently,
\[
 H(A,\theta)=
 \begin{cases}
 A,&A\cot\theta\le\cos^2\theta,\\
 \theta+\sqrt{A\cot\theta(1-A\cot\theta)}-\arccos\sqrt{A\cot\theta},
 &\cos^2\theta<A\cot\theta<1,\\
 \theta,&A\cot\theta\ge1,
 \end{cases}
\]
In the intermediate region,
\[
 \begin{aligned}
 H_{AA}&=-\frac{\cot^2\theta}{2(A\cot\theta)^{3/2}\sqrt{1-A\cot\theta}}<0,\\
 \det D^2H&=\frac{1+\cot^2\theta}{A\cot\theta}
 \bigl((1+\cot^2\theta)A\cot\theta-\cot^2\theta\bigr)>0.
 \end{aligned}
\]
The interface gradients are $(1,0)$ and $(0,1)$, so $H$ is jointly concave.
Include degenerate triangles as needed, using continuity at $(0,0)$. Since
\[
 \sum_{\nu=1}^6A_\nu=\abs P,\qquad
 \sum_{\nu=1}^6\theta_\nu=\pi,
\]
Jensen's inequality gives
\[
 \abs{P\cap B_1}=\sum_{\nu=1}^6 H(A_\nu,\theta_\nu)
 \le6H(\abs P/6,\pi/6)=\abs{H\cap B_1}.
\]
For a nondecreasing radial profile $f$, Stieltjes integration gives
\[
 \int_Pf-\int_Hf
 =\int_0^\infty\bigl(\abs{H\cap B_r}-\abs{P\cap B_r}\bigr)\dd f(r)\ge0.
\]
For nonincreasing $f$, apply this identity to $-f$.
\end{proof}

Use the radial profile $y\mapsto g(\sqrt{\pi a}|y|)$, where
\begin{equation}\label{eq:quintic-profile}
 g(s)=s-\frac5{12}s^3+\frac1{20}s^5\quad(0\le s\le1),
 \qquad g(s)=\frac{19}{30}\quad(s\ge1).
\end{equation}
On $[0,1]$,
\[
 g'=(1-s^2)(1-s^2/4)\ge0,\qquad
 g''=s(s^2-5/2)\le0,\qquad
 (g/s)'=s(s^2/5-5/6)\le0.
\]
Thus
\[
 (g^2)'\ge0,\qquad
 \frac{\dd}{\dd s}\left(g'(s)^2+\frac{g(s)^2}{s^2}\right)\le0,
 \qquad g'(s)^2+\frac{g(s)^2}{s^2}=\frac{361}{900s^2}\quad(s\ge1).
\]
By \eqref{eq:projected} and Lemma~\ref{lem:disk},
\[
 \mu_1(V_\Lambda)\le\pi a
 \frac{\int_H\left(g'(|y|)^2+g(|y|)^2/|y|^2\right)\dd y}
      {\int_Hg(|y|)^2\dd y},
\]
where $H$ now has area $\pi$.

For $0\le s\le1$,
\[
 \begin{aligned}
 \frac gs&=\frac{19}{60}(3-s^2)+\frac{(1-s^2)^2}{20},\\
 g'(s)^2+\frac{g(s)^2}{s^2}-\frac{361}{900}
 &\ge\frac{361}{900}\bigl((1-s^2)+(1-s^2)^2\bigr).
 \end{aligned}
\]
For $s\ge1$,
\[
 \frac{361}{900}-g'(s)^2-\frac{g(s)^2}{s^2}
 =\frac{361}{900}\frac{s^2-1}{s^2}
 \ge\frac{361}{900}\bigl((s^2-1)-(s^2-1)^2\bigr).
\]
Since $|B_1|=|H|$,
\begin{align*}
 &\int_{B_1}\left(g'(|y|)^2+\frac{g(|y|)^2}{|y|^2}\right)\dd y
 -\int_H\left(g'(|y|)^2+\frac{g(|y|)^2}{|y|^2}\right)\dd y\\
 &\quad\ge\frac{361}{900}\left(
 \int_H(3|y|^2-|y|^4)-\int_{B_1}(3|y|^2-|y|^4)\right)\\
 &\quad=\frac{361\pi}{900}\left(\frac{5\pi}{6\sqrt3}-\frac{14\pi^2}{405}-\frac76\right).
\end{align*}
With $\Nzero=\{0,1,2,\ldots\}$, the barycentric coordinates on a facet triangle $\triangle$ give
\[
 |y|^2=\frac{2\pi}{3\sqrt3}(\lambda^2+\lambda\mu+\mu^2),\qquad
 \frac1{|\triangle|}\int_\triangle\lambda^i\mu^j
 =\frac{2i!j!}{(i+j+2)!}\quad(i,j\in\Nzero).
\]
Consequently,
\[
 \frac1\pi\int_H|y|^2=\frac{5\pi}{18\sqrt3},\qquad
 \frac1\pi\int_H|y|^4=\frac{14\pi^2}{405},\qquad
 \int_{B_1}|y|^{2k}=\frac\pi{k+1}\quad(k\in\Nzero).
\]
The disk integrals are
\[
 \int_0^1\left(sg'(s)^2+\frac{g(s)^2}{s}\right)\dd s=\frac{6469}{13500},\qquad
 \int_0^1sg(s)^2\dd s=\frac{1357}{9600}.
\]
Lemma~\ref{lem:disk} gives $\int_Hg^2\ge1357\pi/4800$, and hence
\begin{equation}\label{eq:kappa}
 \mu_1(V_\Lambda)\le a\left(
 \frac{34232\pi}{6785}+\frac{80864\pi^3}{1648755}
 -\frac{14440\pi^2}{12213\sqrt3}\right)<\frac{5317}{500}a.
\end{equation}
For $a=\lambda_1(\Lambda(p))$, Lemma~\ref{lem:lengths} yields
\begin{equation}\label{eq:primalbound}
 F(p)\le\frac{5317}{500}D^{-1/6}
 (\min_{\varnothing\ne S\subsetneq I}W_S)^{1/2}.
\end{equation}
The rational constants also satisfy
\begin{equation}\label{eq:cornerconstants}
 \frac{132379}{100000}\pi^6<T^3,
 \qquad \frac{28}{25}\left(\frac{5317}{500}\right)^6<T^6.
\end{equation}

\section{The spectrum of the BCC cell}\label{sec:spectrum}
Parity decomposition identifies the three axial first modes of the BCC cell; a mixed-problem residual estimate gives \eqref{sp:bcc-enclosure}.

Let
\begin{equation}\label{sp:cell}
 W=[-1,1]^3\cap\{\abs{x_1}+\abs{x_2}+\abs{x_3}\leq 3/2\},
 \qquad \Omega_+=W\cap[0,1]^3.
\end{equation}
The bisectors of $\pm2^{1/3}e_i$ and
$2^{-2/3}(\pm1,\pm1,\pm1)$ give
\[
 V_{\BCC}\subset2^{-2/3}W,\qquad
 \abs{\Omega_+}=\frac12,
 \qquad \abs{V_{\BCC}}=\abs{2^{-2/3}W}=1.
\]
Hence $V_{\BCC}=2^{-2/3}W$. The vertices of $W$ are the signed
permutations of $(1,1/2,0)$; its square and regular hexagonal facets have
common edge length $1/\sqrt2$.

Opposite parities in any coordinate give
\[
 \int_Wuv=0,\qquad \int_W\nabla u\cdot\nabla v=0.
\]
Coordinate reflections therefore split the Neumann spectrum into eight
orthogonal sectors. On $\Omega_+$, odd and even parity in $x_i$ give
Dirichlet and Neumann conditions on $x_i=0$, respectively. Let $\lambda,\lambda_2$ be the first two
eigenvalues with Dirichlet conditions on $x_1=0$ and Neumann conditions elsewhere.

\begin{lemma}\label{sp:separation}
\begin{equation}\label{sp:separation-bounds}
 \lambda\geq\frac{\pi^2}{4},\qquad \lambda_2>\pi^2.
\end{equation}
\end{lemma}

\begin{proof}
For $f\in H^1(0,1)$ with $f(0)=0$, the Dirichlet--Neumann inequality follows from
\[
 \begin{aligned}
 \int_0^1\left(|f'|^2-\frac{\pi^2}{4}|f|^2\right)\dd x
 &=\int_0^1\sin^2\frac{\pi x}{2}
 \left|\left(\frac{f}{\sin(\pi x/2)}\right)'\right|^2\dd x\geq0.
 \end{aligned}
\]
Rescaling on each $x_1$-fibre gives
\[
 \int_{\Omega_+}|\partial_{x_1}u|^2
 \geq\int_{\Omega_+}\frac{\pi^2|u|^2}
 {4\min(1,3/2-x_2-x_3)^2}
 \geq\frac{\pi^2}{4}\int_{\Omega_+}|u|^2.
\]
To bound the second eigenvalue, decompose under $x_2\leftrightarrow x_3$ and
restrict to $x_2\geq x_3$.

In the odd sector, $u|_{x_2=x_3}=0$. Slicing in $x_1$ and $(x_2-x_3)/\sqrt2$, with
$s=x_2+x_3$, gives
\[
 \frac{\int|\nabla u|^2}{\int|u|^2}
 \geq\inf_{0<s<3/2}\pi^2\left\{
 \frac{1}{4\min(1,3/2-s)^2}
 +\frac{1}{2\min(s,2-s)^2}\right\}\geq\frac54\pi^2.
\]
Indeed, H\"older's inequality yields
\[
 \frac{1}{4\min(1,3/2-s)^2}
 +\frac{1}{2\min(s,2-s)^2}
 \geq
 \begin{cases}
 9/4,&0<s\leq1/2,\\
 (1+2^{1/3})^3/9>5/4,&1/2<s<1,\\
 3/2,&1\leq s<3/2.
 \end{cases}
\]

In the even sector, insert a Neumann cut at $x_1=1/2$ and put
\[
 \Omega_L=\Omega_+\cap\{x_2\geq x_3,\ x_1<1/2\},\qquad
 \Omega_R=\Omega_+\cap\{x_2\geq x_3,\ x_1>1/2\}.
\]
On $\Omega_L$, Dirichlet--Neumann fibres have length at most $1/2$, with
strict inequality on a set of positive transverse measure. Positivity of
the ground state gives $\lambda_1(\Omega_L)>\pi^2$. On $\Omega_R$, set
\[
 (x_2,x_3)=(3/2-x_1)(Y,Z),\qquad
 T_0=\{0\leq Z\leq Y,\ Y+Z\leq1\},\qquad |T_0|=\frac14.
\]
The triangle and its reflection across $Z=0$ form a square of side $1/\sqrt2$. Thus
\[
 dV=(3/2-x_1)^2\,dx_1\,dY\,dZ,\qquad
 \int_{T_0}\left|f-4\int_{T_0}f\right|^2
 \leq\frac1{2\pi^2}\int_{T_0}|\nabla f|^2.
\]

For the first nonconstant weighted Neumann mode,
\[
 -\big((3/2-x_1)^2f'\big)'=\nu(3/2-x_1)^2f,\qquad
 f'(1/2)=f'(1)=0,
\]
and hence
\[
 \begin{aligned}
 \nu\int_{1/2}^1(3/2-x_1)^2(f')^2
 &=\int_{1/2}^1(3/2-x_1)^2(f'')^2+2\int_{1/2}^1(f')^2\\
 &\geq\int_{1/2}^1\bigl(((3/2-x_1)f')'\bigr)^2
 \geq4\pi^2\int_{1/2}^1(3/2-x_1)^2(f')^2.
 \end{aligned}
\]
Integration by parts gives
\[
 \int_{1/2}^1\bigl(((3/2-x_1)f')'\bigr)^2
 =\int_{1/2}^1(3/2-x_1)^2(f'')^2.
\]
Set
\[
 m(x_1)=4\int_0^1\int_0^{\min(Y,1-Y)}
 u\bigl(x_1,(3/2-x_1)Y,(3/2-x_1)Z\bigr)\dd Z\dd Y.
\]
Then
\[
 \begin{aligned}
 4\int_{T_0}(Y^2+Z^2)\dd Y\dd Z&=\frac13,\\
 m'(x_1)&=4\int_{T_0}(u_{x_1}-Yu_{x_2}-Zu_{x_3})\dd Y\dd Z,\\
 |m'(x_1)|^2&\leq8\int_{T_0}u_{x_1}^2\dd Y\dd Z
 +\frac83\int_{T_0}(u_{x_2}^2+u_{x_3}^2)\dd Y\dd Z,
 \end{aligned}
\]
where the derivatives are evaluated at $\bigl(x_1,(3/2-x_1)Y,(3/2-x_1)Z\bigr)$.
Since
\[
 \begin{gathered}
 |\Omega_R|=\frac7{96},\qquad
 \frac{96}{7}\int_{\Omega_R}u
 =\frac{24}{7}\int_{1/2}^1(3/2-x_1)^2m(x_1)\dd x_1,\\
 \frac12\leq3/2-x_1\leq1,
 \end{gathered}
\]
the transverse and weighted inequalities give
\begin{align*}
 \int_{\Omega_R}|u-m|^2
 &\leq\frac1{2\pi^2}\int_{\Omega_R}(3/2-x_1)^2(u_{x_2}^2+u_{x_3}^2),\\
 \int_{\Omega_R}\left|m-\frac{96}{7}\int_{\Omega_R}u\right|^2
 &\leq\frac1{16\pi^2}\int_{1/2}^1(3/2-x_1)^2|m'|^2\dd x_1\\
 &\leq\frac1{2\pi^2}\int_{\Omega_R}u_{x_1}^2
 +\frac1{6\pi^2}\int_{\Omega_R}(u_{x_2}^2+u_{x_3}^2).
\end{align*}
Thus
\[
 \begin{aligned}
 \int_{\Omega_R}\left|u-\frac{96}{7}\int_{\Omega_R}u\right|^2
 &=\int_{\Omega_R}|u-m|^2
 +\int_{\Omega_R}\left|m-\frac{96}{7}\int_{\Omega_R}u\right|^2\\
 &\leq\frac2{3\pi^2}\int_{\Omega_R}|\nabla u|^2.
 \end{aligned}
\]
Neumann bracketing therefore gives
\[
 \begin{gathered}
 \lambda_1(\Omega_R)=0,\qquad
 \lambda_2(\Omega_R)\geq\frac32\pi^2,\\
 \lambda_2^{\mathrm{even}}
 \geq\min\{\lambda_1(\Omega_L),\lambda_2(\Omega_R)\}>\pi^2.
 \end{gathered}
\]
The odd-sector bound completes \eqref{sp:separation-bounds}.
\end{proof}

For complementary bounds using $H(\operatorname{div})$ fluxes, see
\v{S}ebestov\'a and Vejchodsk\'y
\cite{SebestovaVejchodsky} and Vejchodsk\'y \cite{Vejchodsky}.
Assume
\[
 \begin{gathered}
 v\in H^1(\Omega_+),\qquad \norm v_{L^2(\Omega_+)}=1,
 \qquad v|_{x_1=0}=0,\\
 q\in H(\operatorname{div};\Omega_+),\qquad
 (q\cdot n)|_{\partial\Omega_+\setminus\{x_1=0\}}=0.
 \end{gathered}
\]
\begin{lemma}\label{sp:residual}
If $\norm{\nabla v}_2^2<\gamma\leq\lambda_2$, then
\begin{equation}\label{sp:residual-bound}
 \begin{aligned}
 0\leq\norm{\nabla v}_2^2-\lambda
 &\leq\frac{\gamma}{\gamma-\norm{\nabla v}_2^2}
 \Bigl(\frac{4001}{4000}\norm{\nabla v-q}_2^2\\
 &\hspace{65pt}+\frac{16004}{9}
 \norm{\operatorname{div}q+\norm{\nabla v}_2^2v}_2^2\Bigr).
 \end{aligned}
\end{equation}
\end{lemma}

\begin{proof}
For every mixed test function $w$, integration by parts and
Lemma~\ref{sp:separation} give
\[
 \begin{aligned}
 &\left|\int_{\Omega_+}
   (\nabla v\cdot\nabla w-\norm{\nabla v}_2^2vw)\right|\\
 &\quad=\left|\int_{\Omega_+}
  \bigl((\nabla v-q)\cdot\nabla w
       -(\operatorname{div}q+\norm{\nabla v}_2^2v)w\bigr)\right|\\
 &\quad\leq\Bigl(\norm{\nabla v-q}_2
       +\frac23\norm{\operatorname{div}q+\norm{\nabla v}_2^2v}_2\Bigr)
       \norm{\nabla w}_2.
 \end{aligned}
\]
In an orthonormal mixed eigenbasis,
\[
 \begin{gathered}
 v=\sum_{j=1}^{\infty}v_j\phi_j,\qquad
 \sum_{j=1}^{\infty}|v_j|^2=1,\qquad
 \norm{\nabla v}_2^2=\sum_{j=1}^{\infty}\lambda_j|v_j|^2,\\
 (\gamma-\norm{\nabla v}_2^2)(1-|v_1|^2)
 \leq(\norm{\nabla v}_2^2-\lambda)|v_1|^2,\qquad
 |v_1|^2\geq\frac{\gamma-\norm{\nabla v}_2^2}{\gamma-\lambda}.
 \end{gathered}
\]
Duality and Cauchy--Schwarz yield
\begin{align*}
 &\frac{4001}{4000}\norm{\nabla v-q}_2^2
   +\frac{16004}{9}\norm{\operatorname{div}q+\norm{\nabla v}_2^2v}_2^2\\
 &\quad\geq\sum_{j=1}^{\infty}
       \frac{(\lambda_j-\norm{\nabla v}_2^2)^2}{\lambda_j}|v_j|^2\\
 &\quad\geq(\norm{\nabla v}_2^2-\lambda)|v_1|^2
       \left(\frac{\norm{\nabla v}_2^2-\lambda}{\lambda}
             +1-\frac{\norm{\nabla v}_2^2}{\gamma}\right)\\
 &\quad\geq\frac{\norm{\nabla v}_2^2-\lambda}{\lambda}
       \frac{\norm{\nabla v}_2^2(\gamma-\norm{\nabla v}_2^2)}{\gamma}\\
 &\quad\geq(\norm{\nabla v}_2^2-\lambda)
       \frac{\gamma-\norm{\nabla v}_2^2}{\gamma}.
\end{align*}
This proves \eqref{sp:residual-bound}.
\end{proof}

Set $b=(1/2,1/2,1/2)$ and $v_{ij}=e_i+e_j/2$. Coning the
non-coordinate faces of $\Omega_+$ to $0$ gives the nine tetrahedra of
$\mathcal T_+$:
\begin{equation}\label{sp:tetrahedra}
 \begin{aligned}
 &(0,e_i,v_{ij},v_{ik}),\quad (0,b,v_{ij},v_{ik})
 &&\bigl(i=1,2,3,\ j<k,\ \{i,j,k\}=\{1,2,3\}\bigr),\\
 &(0,b,v_{ij},v_{ji})&& (1\leq i<j\leq3).
 \end{aligned}
\end{equation}
Order the first six by increasing $i$, with the square cone first for each $i$, and the last
three lexicographically in $(i,j)$.
For barycentric coordinates $\zeta$ on $P$, use
\[
 \alpha\in\Nzero^4,\qquad
 |\alpha|=\sum_{i=0}^3\alpha_i,\qquad
 \alpha!=\prod_{i=0}^3\alpha_i!.
\]
Then
\begin{equation}\label{sp:barycentric-integral}
 \int_P\zeta^\alpha
 =6\abs P\frac{\prod_{i=0}^3\alpha_i!}{(3+\abs\alpha)!},\qquad \alpha\in\Nzero^4.
\end{equation}
For the trial and flux construction, use Bernstein coefficients:
\[
 f|_P=\sum_{\substack{\alpha\in\Nzero^4\\|\alpha|=m}}
 f_{P,\alpha}\frac{m!}{\alpha!}\zeta^\alpha.
\]

For $f,h$ of degrees $m,n$, \eqref{sp:barycentric-integral} gives
\begin{equation}\label{sp:bernstein-product}
 \int_Pfh=\frac{6\abs{P}m!n!}{(m+n+3)!}
 \sum_{\substack{\alpha\in\Nzero^4\\\abs\alpha=m}}
 \sum_{\substack{\beta\in\Nzero^4\\\abs\beta=n}}
 f_\alpha h_\beta
 \prod_{i=0}^3\binom{\alpha_i+\beta_i}{\alpha_i}.
\end{equation}
In particular,
\[
 (f,h)_{\Omega_+}=\sum_{P\in\mathcal T_+}\int_Pfh.
\]
\begin{numlemma}\label{num:flux}
\[
 \begin{gathered}
 \exists v\in H^1(\Omega_+)\cap C(\overline{\Omega_+}),\qquad
 \exists q\in H(\operatorname{div};\Omega_+),\\
 \norm v_2=1,\qquad v|_{x_1=0}=0,\qquad
 (q\cdot n)|_{\partial\Omega_+\setminus\{x_1=0\}}=0,\\
 v|_P,\ q_i|_P\in\R[x_1,x_2,x_3],\qquad
 \deg(v|_P)\le9,\quad\deg(q_i|_P)\le10
 \quad(P\in\mathcal T_+,\ 1\le i\le3),
 \end{gathered}
\]
\begin{equation}\label{sp:trial-bounds}
 \begin{aligned}
 \frac{43008469}{10^7}<\norm{\nabla v}_2^2
 &<\frac{43008470}{10^7},\\
 \norm{\nabla v-q}_2^2&<\frac{9002}{10^9},\\
 \left\lVert\operatorname{div}q+\norm{\nabla v}_2^2v\right\rVert_2^2
 &<\frac{23}{10^{13}}.
 \end{aligned}
\end{equation}
\end{numlemma}

The Appendix gives the construction and finite rational evaluations.
Taking $\gamma=49/5$ in Lemma~\ref{sp:residual} and using
\eqref{sp:separation-bounds} and \eqref{sp:trial-bounds} gives
\begin{equation}\label{sp:mixed-enclosure}
 \frac{43008308}{10^7}<\lambda<\frac{43008470}{10^7}.
\end{equation}

\begin{proposition}\label{sp:axial}
\[
 \begin{gathered}
 \mu_1(W)=\lambda,\qquad
 \ker(-\Delta_W^N-\lambda)=\operatorname{span}\{u_1,u_2,u_3\},\\
 \int_Wu_iu_j=\delta_{ij},\qquad
 u_i(x-2x_je_j)=(-1)^{\delta_{ij}}u_i(x),\\
 u_i(x)>0\quad(x\in\operatorname{int}W,\ x_i>0),
 \qquad 1\le i,j\le3.
 \end{gathered}
\]
\begin{equation}\label{sp:bcc-enclosure}
 \mu_1(V_{\BCC})=2^{4/3}\lambda,
 \qquad 10.837414<\mu_1(V_{\BCC})<10.837456.
\end{equation}
\end{proposition}

\begin{proof}
By \eqref{sp:mixed-enclosure}, $\lambda<\pi^2/2$. For $u$ odd in
at least two coordinates, slicing gives
\[
 \int_W|\nabla u|^2\geq\frac{\pi^2}{2}\int_W|u|^2.
\]

Let $u$ be a mean-zero all-even eigenfunction with $\nu<\pi^2/2$.
For a coordinate transposition $\tau$,
\[
 \Delta(u\circ\tau)=(\Delta u)\circ\tau,\qquad
 -\Delta(u-u\circ\tau)=\nu(u-u\circ\tau).
\]
The difference is odd in the corresponding difference coordinate, whose
half-fibres have length at most $1/\sqrt2$. Hence
\[
 u-u\circ\tau\neq0
 \quad\Longrightarrow\quad
 \nu=\frac{\int_{\Omega_+}|\nabla(u-u\circ\tau)|^2}
 {\int_{\Omega_+}|u-u\circ\tau|^2}\geq\frac{\pi^2}{2}.
\]
Thus $u=u\circ\tau$ for every coordinate transposition, and
\[
 \begin{gathered}
 \operatorname{Vert}\bigl(\Omega_+\cap\{x_1\geq x_2\geq x_3\}\bigr)\\
 =\{0,\ (1,0,0),\ (3/4,3/4,0),\ (1/2,1/2,1/2),
       \ (1,1/2,0),\ (1,1/4,1/4)\}.
 \end{gathered}
\]
The Payne--Weinberger inequality \cite{PayneWeinberger} gives
\[
 \begin{gathered}
 \operatorname{diam}\bigl(\Omega_+\cap\{x_1\geq x_2\geq x_3\}\bigr)
 =\frac{\sqrt5}{2},\\
 \int_{\Omega_+\cap\{x_1\geq x_2\geq x_3\}}u
 =\frac16\int_{\Omega_+}u=0,\\
 \nu=\frac{\displaystyle\int_{\Omega_+\cap\{x_1\geq x_2\geq x_3\}}|\nabla u|^2}
 {\displaystyle\int_{\Omega_+\cap\{x_1\geq x_2\geq x_3\}}|u|^2}
 \geq\frac{4\pi^2}{5}>\frac{\pi^2}{2}.
 \end{gathered}
\]
Simplicity and positivity of the mixed ground state give
\[
 \mu_1(W)=\lambda,\qquad
 \ker(-\Delta_W^N-\lambda)=\operatorname{span}\{u_1,u_2,u_3\},
 \qquad \mu_1(V_{\BCC})=2^{4/3}\lambda.
\]
Scaling \eqref{sp:mixed-enclosure} gives \eqref{sp:bcc-enclosure}
and $\mu_1(V_{\BCC})>T$.
\end{proof}

The ten-place numerical value in \eqref{eq:bcc-value} was computed using
conforming high-order approximations on geometrically graded meshes and checked under
polynomial enrichment, mesh refinement, and increased quadrature order.
The proof uses the rigorous enclosure \eqref{sp:mixed-enclosure}.

\section{Strict comparison near BCC}\label{sec:local}
We prove strict local maximality of BCC by transporting the axial trial space and controlling its linear splitting and higher-order terms.

\begin{proposition}\label{prop:local}
\[
 0<\lVert p-p_*\rVert_2\leq\frac{73}{6000}
 \quad\Longrightarrow\quad F(p)<F(p_*).
\]
\end{proposition}
For the cell $W$ of volume $4$ in Section~\ref{sec:spectrum}, fix
\[
 (b_0,b_1,b_2,b_3)=((1,1,1),(1,-1,-1),(-1,1,-1),(-1,-1,1)).
\]
The superbase metric is $I+K$, where
\begin{equation}\label{lo:strain}
 \begin{gathered}
 K=\frac1{16}\sum_{0\le i<j\le3}(6p_{ij}-1)(b_i-b_j)(b_i-b_j)^T,
 \qquad\tr K=0,\\
 36\lVert p-p_*\rVert_2^2
 =8\lVert\diag K\rVert_F^2+4\lVert\off K\rVert_F^2.
 \end{gathered}
\end{equation}
\[
 \sum_{i=0}^3b_i b_i^T=4I,\qquad b_i^TKb_j=1-6p_{ij}\quad(i\ne j),
\]
and hence
\[
 \begin{gathered}
 16\lVert K\rVert_F^2=\sum_{i=0}^3\sum_{j=0}^3(b_i^TKb_j)^2
 =72\lVert p-p_*\rVert_2^2+8\lVert\off K\rVert_F^2,\\
 \frac3{\sqrt2}\lVert p-p_*\rVert_2
 \leq\lVert K\rVert_F\leq3\lVert p-p_*\rVert_2.
 \end{gathered}
\]
In the hyperplane $\sum_{i=0}^3y_i=0$, use $x=\frac14\sum_{i=0}^3y_i b_i$. A vertex indexed by a strict order $\sigma$ is
\begin{equation}\label{lo:vertices}
 x_\sigma(p)=\frac34\sum_{0\le i<j\le3}p_{ij}
 \sgn(\sigma_i-\sigma_j)(b_i-b_j).
\end{equation}
Since $6p_{ij}\geq1-6\lVert p-p_*\rVert_2>0$, the face incidences persist. Triangulate each face from its vertex average and cone the resulting triangles to $0$; matching the vertices and face averages of the resulting $72$ tetrahedra defines an odd piecewise linear map with matrix $I+E$ on each tetrahedron.

For the fan \eqref{sp:tetrahedra} of Section~\ref{sec:spectrum}, the vertex displacements satisfy
\begin{align}\label{lo:displacements}
 \delta v_{ij}&=K_{ii}e_i+(K_{ij}-K_{ii})e_j+(K_{ik}+K_{jk})e_k,
 &\delta e_i&=Ke_i,&\delta b&=Kb.
\end{align}
Thus, on $T=(0,c,v,w)$,
\begin{equation}\label{lo:transportmatrix}
 E_T(K)=[\delta c,\delta v,\delta w][c,v,w]^{-1}.
\end{equation}
For a coordinate reflection $S$, the matrix on $ST$ is $SE_T(SKS)S$. We claim
\begin{equation}\label{lo:Ebounds}
 \begin{gathered}
 \lVert E\rVert_F\leq6\lVert p-p_*\rVert_2,\\
 |\tr E_{\rm sq}|^2\leq27\lVert p-p_*\rVert_2^2,\qquad
 |\tr E_{\rm hex}|^2\leq12\lVert p-p_*\rVert_2^2.
 \end{gathered}
\end{equation}
The active vectors at $v_{ij}$ are $2e_i$ and $e_i+e_j\pm e_k$. Differentiating their bisector equations $v\cdot x=v^T(I+K)v/2$ gives
\[
 2e_i\cdot\delta v_{ij}=2K_{ii},\qquad
 (e_i+e_j\pm e_k)\cdot\delta v_{ij}=K_{ij}\pm(K_{ik}+K_{jk}).
\]
These equations give \eqref{lo:displacements}; face averaging gives $\delta e_i=Ke_i$, $\delta b=Kb$. Since $\tr K=0$,
\[
 36\lVert p-p_*\rVert_2^2=12K_{ii}^2+4(K_{jj}-K_{kk})^2+8(K_{ij}^2+K_{ik}^2+K_{jk}^2).
\]
The square cones give
\[
 E_{\rm sq}=\begin{pmatrix}
 K_{ii}&0&0\\K_{ij}&-2K_{ii}&2K_{jk}\\K_{ik}&2K_{jk}&-2K_{ii}
 \end{pmatrix}.
\]
In the basis $e_i,(e_j+e_k)/\sqrt2,(e_j-e_k)/\sqrt2$, the two hexagon edge pairs $(v_{ij},v_{ik})$ and $(v_{jk},v_{kj})$ give
\[
 E_{\rm hex}=\begin{pmatrix}
 K_{ii}-(K_{ij}+K_{ik})/3&4(K_{ij}+K_{ik})/(3\sqrt2)&0\\
 (-K_{ii}+K_{ij}+K_{ik}+2K_{jk}/3)/\sqrt2&2K_{jk}/3&0\\
 (K_{ij}-K_{ik}-(K_{jj}-K_{kk})/3)/\sqrt2&2(K_{jj}-K_{kk})/3&-2(K_{ii}+K_{jk})
 \end{pmatrix}.
\]
Consequently,
\begin{align*}
 36\lVert p-p_*\rVert_2^2-\lVert E_{\rm sq}\rVert_F^2
 &=3K_{ii}^2+4(K_{jj}-K_{kk})^2+7K_{ij}^2+7K_{ik}^2\geq0,\\
 36\lVert p-p_*\rVert_2^2-\lVert E_{\rm hex}\rVert_F^2
 &=\begin{pmatrix}K_{ii}&K_{ij}+K_{ik}&K_{jk}\end{pmatrix}
 \begin{pmatrix}13/2&5/6&-11/3\\5/6&5/2&-1/3\\-11/3&-1/3&10/3\end{pmatrix}
 \begin{pmatrix}K_{ii}\\K_{ij}+K_{ik}\\K_{jk}\end{pmatrix}\\
 &\quad+\begin{pmatrix}K_{jj}-K_{kk}&K_{ij}-K_{ik}\end{pmatrix}
 \begin{pmatrix}7/2&1/6\\1/6&7/2\end{pmatrix}
 \begin{pmatrix}K_{jj}-K_{kk}\\K_{ij}-K_{ik}\end{pmatrix}\geq0.
\end{align*}
The two coefficient blocks are positive definite. Also,
\begin{align*}
 |\tr E_{\rm sq}|^2&=9K_{ii}^2\leq27\lVert p-p_*\rVert_2^2,\\
 |\tr E_{\rm hex}|^2
 &=\left(K_{ii}+\frac{K_{ij}+K_{ik}+4K_{jk}}3\right)^2
 \leq12\lVert p-p_*\rVert_2^2.
\end{align*}
This proves \eqref{lo:Ebounds}.

The map $\det(I+K)^{-1/6}(I+K)^{-1/2}$ gives the physical volume-four cell. Its constant Jacobian cancels from the Rayleigh quotient, and its inverse metric is $\det(I+K)^{1/3}(I+K)$, with
\[
 \det(I+K)^{1/3}
 \leq\frac{\tr(I+K)}3=1,\qquad
 \lVert E\rVert_{\rm op}\leq6\lVert p-p_*\rVert_2<1.
\]
Put
\begin{equation}\label{lo:JC}
 J=\det(I+E)>0,\qquad
 C=J(I+E)^{-1}(I+K)(I+E)^{-T}.
\end{equation}
Define
\[
 \mathcal B_\theta(f,g;A,b)=\int_W(\nabla f^TA\nabla g-\theta bfg).
\]
For $0\ne f\in H^1(W)$ with $f(-x)=-f(x)$,
\[
 J(-x)=J(x),\qquad \int_WJf=0,
\]
so the Rayleigh principle gives
\begin{equation}\label{lo:Rayleigh}
 2^{-4/3}F(p)-\theta
 \leq\frac{\mathcal B_\theta(f,f;C,J)}{\int_WJf^2},
 \qquad \int_WJf^2>0.
\end{equation}

On $\Omega_+$ define
\begin{equation}\label{lo:polynomial}
 \begin{aligned}
 p_1={}&x_1\sum_{\substack{i,j,k\in\Nzero\\i+j+k\le7,\ j\le k}}c_{ijk}
 \sum_{a=0}^i\sum_{b=0}^j\sum_{c=0}^k
 (-1)^{i+j+k-a-b-c}\\
 &\quad\times\binom ia\binom{i+a}a\binom jb\binom{j+b}b
 \binom kc\binom{k+c}c\,
 x_1^a\bigl(x_2^bx_3^c+\1_{j<k}x_2^cx_3^b\bigr).
 \end{aligned}
\end{equation}
Read the $70$ coefficients row by row, from left to right and top to bottom,
in the nested order
\[
 i=0,\ldots,7,\qquad
 j=0,\ldots,\left\lfloor\frac{7-i}{2}\right\rfloor,\qquad
 k=j,\ldots,7-i-j.
\]
\begingroup\small
\[
\begin{array}{rrrrr}
229900428&-139871184&-52930218&-14629008&-13261341\\
-5686815&-889274&73926&-83267960&-41119660\\
-26184507&-21182213&-8427964&-1388016&-11380873\\
-6083225&-7172094&-2162237&282052&-402975\\
-775093&78837850&21184188&-10453204&-10599016\\
-5677196&-1112077&116720818&67372472&14741791\\
-3555916&-3386892&47315301&14147076&2266381\\
2353273&-3445608&11601610&-2439555&-3913125\\
-1623965&-757041&-10857375&-10285059&-624273\\
1512055&-7195789&-1214356&9282897&10229380\\
4231107&1067499&337892&854262&-3570128\\
-1636479&-674516&8206128&11338133&3833423\\
487474&13246110&3072687&3185632&4391020\\
1201058&4365604&340142&95828&-177693
\end{array}
\]
\endgroup
Thus $c_{000}=229900428$, $c_{005}=-5686815$, and $c_{700}=-177693$.
Obtain $p_i$ by interchanging coordinates $1$ and $i$. Set
\begin{equation}\label{lo:normalization}
 \begin{gathered}
 \int_{\Omega_+}p_1^2
 =\frac{6755696638648162400908255727}{675569664000},\\
 P_i(x)=\frac{\sgn(x_i)}{\sqrt{8\int_{\Omega_+}p_1^2}}p_i(|x_1|,|x_2|,|x_3|),\qquad x\in W.
 \end{gathered}
\end{equation}

For $\alpha\in\Nzero^3$, write $|\alpha|=\alpha_1+\alpha_2+\alpha_3$ and $\alpha!=\alpha_1!\alpha_2!\alpha_3!$. Barycentric integration on $T=\operatorname{conv}(w_0,\ldots,w_3)\in\mathcal T_+$ gives
\begin{equation}\label{lo:moments}
 \int_Tx^\alpha\,dx
 =\frac{6|T|\alpha!}{(|\alpha|+3)!}
 [z^\alpha]\prod_{i=0}^3(1-w_i\cdot z)^{-1}.
\end{equation}
For the doubled vertices, the identities
\[
 \sum_{\alpha\in\Nzero^3}g_\alpha^{(0)}z^\alpha=1,\qquad
 (1-2w_i\cdot z)\sum_{\alpha\in\Nzero^3}g_\alpha^{(i+1)}z^\alpha
 =\sum_{\alpha\in\Nzero^3}g_\alpha^{(i)}z^\alpha
\]
give $g_\alpha^{(0)}=\boldsymbol1_{\alpha=0}$ and, for $i=0,1,2,3$,
\begin{equation}\label{lo:momentrecurrence}
 g_\alpha^{(i+1)}=g_\alpha^{(i)}+
 \sum_{j=1}^3(2w_i)_j g_{\alpha-e_j}^{(i+1)},
 \qquad g_\beta^{(i+1)}=0\quad\text{if }\min_{1\le j\le3}\beta_j<0.
\end{equation}
Evaluation in increasing $|\alpha|$ through $16$ gives
\begin{gather*}
 [z^\alpha]\prod_{i=0}^3(1-w_i\cdot z)^{-1}=2^{-|\alpha|}g_\alpha^{(4)},\\
 8(6|T|)\in\Z,\qquad
 2^{19}19!\int_Tx^\alpha\,dx\in\Z\qquad(|\alpha|\leq16).
\end{gather*}
Summing these moments against the coefficients of $p_1^2$ and $|\nabla p_1|^2$ yields \eqref{lo:normalization} and \eqref{lo:orthogonality}.

Parity and coordinate permutation give
\begin{equation}\label{lo:orthogonality}
 \begin{gathered}
 \int_WP_iP_j=\delta_{ij},\qquad
 \int_W\nabla P_i\cdot\nabla P_j=r\delta_{ij},\\
 r=\frac{377741884455339525044322247762}{87824056302426111211807324451}
 <\frac{430113}{100000}.
 \end{gathered}
\end{equation}
By \eqref{sp:mixed-enclosure}, $4300827/10^6<\lambda$.
With unit octant normalization and a common sign choice for the eigenfunctions, write
\[
 P=\alpha u+w,\qquad (w,u)=0,\qquad
 \alpha\geq0,\qquad \alpha^2+\lVert w\rVert_2^2=1.
\]
Since $\lambda_2>\pi^2>49/5$, the mixed spectral gap yields
\begin{align*}
 r-\lambda&=\lVert\nabla w\rVert_2^2-\lambda\lVert w\rVert_2^2
 \geq\frac{549887}{980000}\lVert\nabla w\rVert_2^2,\\
 \lVert\nabla w\rVert_2^2&<\frac1{1800},\qquad
 1-\alpha\leq\lVert w\rVert_2^2<\frac1{17640},\\
 \lVert\nabla(P-u)\rVert_2^2&=\lVert\nabla w\rVert_2^2+\lambda(1-\alpha)^2
 <\left(\frac3{125}\right)^2,\\
 \lVert P-u\rVert_2^2&=\lVert w\rVert_2^2+(1-\alpha)^2
 <\left(\frac1{125}\right)^2.
\end{align*}
For $u_\xi=\sum_{i=1}^3\xi_i u_i$, $P_\xi=\sum_{i=1}^3\xi_iP_i$, $|\xi|=1$, orthogonality of the parity sectors gives
\begin{equation}\label{lo:errors}
 \lVert\nabla(u_\xi-P_\xi)\rVert_2<\frac3{125},
 \qquad \lVert u_\xi-P_\xi\rVert_2<\frac1{125}.
\end{equation}

In dimension three,
\[
 J=1+\tr E+\frac{(\tr E)^2-\tr(E^2)}2+\det E.
\]
For $M_0=I$ and $M_m=0$ when $m<0$, define
\begin{align}\label{lo:homogeneous}
 M_m&=\sum_{a=0}^{m}(-E)^a(-E^T)^{m-a}
       +\sum_{a=0}^{m-1}(-E)^aK(-E^T)^{m-1-a}\qquad(m\geq1),\\
 C_m&=M_m+(\tr E)M_{m-1}+\frac{(\tr E)^2-\tr(E^2)}2M_{m-2}+(\det E)M_{m-3}.\notag
\end{align}
Then $C=\sum_{m=0}^{\infty}C_m$. For $p\ne p_*$,
\eqref{lo:Ebounds} and \eqref{lo:strain} give
\begin{align}\label{lo:homogeneousbounds}
 |\tr E|&<\frac{26}{5}\lVert p-p_*\rVert_2,\\
 \frac12|(\tr E)^2-\tr(E^2)|
 &\leq\frac12\bigl((\tr E)^2+\lVert E\rVert_F^2\bigr)
 \leq\frac{63}{2}\lVert p-p_*\rVert_2^2,\notag\\
 |\det E|&\leq\left(\frac{\lVert E\rVert_F}{\sqrt3}\right)^3
 <\frac{216}{5}\lVert p-p_*\rVert_2^3,\notag\\
 \lVert M_m\rVert_{\rm op}
 &\leq(m+1)(6\lVert p-p_*\rVert_2)^m
 +m\lVert K\rVert_{\rm op}(6\lVert p-p_*\rVert_2)^{m-1}\notag\\
 &\leq(1+3m/2)(6\lVert p-p_*\rVert_2)^m\qquad(m\geq1),\notag\\
 \lVert C_2\rVert_{\rm op}&<261\lVert p-p_*\rVert_2^2,\notag\\
 \lVert C_3\rVert_{\rm op}&<2484\lVert p-p_*\rVert_2^3.\notag
\end{align}
The first homogeneous term is
\[
 C_1=K-E-E^T+(\tr E)I.
\]
Let
\[
 (L_\theta(K/2))_{ij}=\mathcal B_\theta(f_i,f_j;C_1,\tr E),\qquad
 (\theta,f_i)=
 \begin{cases}(\lambda,u_i),\\(r,P_i).\end{cases}
\]
For the eigenfunctions, this is the first-order compression to their eigenspace \cite{Kato}. For every signed permutation $S$,
\[
 L_\theta(SKS^T/2)=SL_\theta(K/2)S^T.
\]
Coordinate reflections and the exchange $1\leftrightarrow2$ give
\[
 \begin{aligned}
 L_\theta(\diag(1,-1,0))
 &=\bigl(L_\theta(\diag(1,-1,0))\bigr)_{11}\diag(1,-1,0),\\
 L_\theta(e_1e_2^T+e_2e_1^T)
 &=\bigl(L_\theta(e_1e_2^T+e_2e_1^T)\bigr)_{12}(e_1e_2^T+e_2e_1^T).
 \end{aligned}
\]
Since $L_\theta$ is linear and $\tr K=0$,
\[
 K=K_{11}\diag(1,0,-1)+K_{22}\diag(0,1,-1)
 +\sum_{1\leq i<j\leq3}K_{ij}(e_ie_j^T+e_je_i^T).
\]
Permutation covariance therefore yields
\begin{equation}\label{lo:splitting}
 \begin{aligned}
 L_\theta(K/2)
 &=\frac12\bigl(L_\theta(\diag(1,-1,0))\bigr)_{11}\diag K\\
 &\quad+\frac12\bigl(L_\theta(e_1e_2^T+e_2e_1^T)\bigr)_{12}\off K.
 \end{aligned}
\end{equation}
For $m\in\{2,3\}$, put
\[
 (B_m(K))_{ij}
 =\mathcal B_r\!\left(P_i,P_j;C_m,[z^m]\det(I+zE)\right).
\]

For $A=A^T$ and pairs $(P,u)$, $(Q,v)$ satisfying \eqref{lo:errors},
expansion gives
\begin{align*}
 &\mathcal B_\lambda(u,v;A,b)-\mathcal B_r(P,Q;A,b)\\
 &\quad=\mathcal B_r(u-P,Q;A,b)+\mathcal B_r(P,v-Q;A,b)\\
 &\qquad+\mathcal B_r(u-P,v-Q;A,b)+(r-\lambda)\int_Wbuv
\end{align*}
and hence, by Cauchy--Schwarz,
\begin{align}\label{lo:transfer}
 &\bigl|\mathcal B_\lambda(u,v;A,b)-\mathcal B_r(P,Q;A,b)\bigr|\\
 &\quad\leq\frac6{125}\max\{\lVert A\nabla P\rVert_2,
                                      \lVert A\nabla Q\rVert_2\}
       +\frac9{15625}\lVert A\rVert_{\infty,\rm op}\notag\\
 &\qquad+\frac{430113}{100000}
       \left(\frac2{125}\max\{\lVert bP\rVert_2,\lVert bQ\rVert_2\}
             +\frac{\lVert b\rVert_\infty}{15625}\right)
       +\frac{303}{10^6}\left|\int_Wbuv\right|.\notag
\end{align}

\begin{numlemma}[BCC trial moments]\label{num:local}
For \eqref{lo:polynomial}--\eqref{lo:normalization}, componentwise,
\begin{equation}\label{lo:numericalcoefficients}
 \begin{pmatrix}
 (L_r(\diag(1,-1,0)))_{11}
   &(L_r(e_1e_2^T+e_2e_1^T))_{12}\\
 (L_\lambda(\diag(1,-1,0)))_{11}
   &(L_\lambda(e_1e_2^T+e_2e_1^T))_{12}
 \end{pmatrix}
 <\begin{pmatrix}-97/50&-263/50\\-7/5&-24/5\end{pmatrix}.
\end{equation}
For $p\ne p_*$,
\begin{equation}\label{lo:highercoefficients}
 \lambda_{\max}(B_2(K))<45\lVert p-p_*\rVert_2^2,
 \qquad \lVert B_3(K)\rVert_F<\frac{2376}{5}\lVert p-p_*\rVert_2^3.
\end{equation}
\end{numlemma}

\begin{proof}[Proof of Numerical Lemma~\ref{num:local}]
Use
\[
 (K,A,b)=
 \begin{cases}
 (2\diag(1,-1,0),C_1+2I,\tr E+2),&a=1,\\
 (2(e_1e_2^T+e_2e_1^T),C_1,\tr E),&a=2.
 \end{cases}
\]
The diagonal shift is valid because $\mathcal B_\theta(f_1,f_1;I,1)=0$. Under coordinate reflections, the diagonal pairing is invariant and the two $x_1x_2$ characters in the shear pairing cancel. Equations~\eqref{lo:transportmatrix} and \eqref{lo:moments} give, componentwise,
\begin{equation}\label{lo:firstorderintegers}
 \begin{gathered}
 \begin{pmatrix}-1941429\\43700852\\8980830\\307200\end{pmatrix}
 <\frac{10^6}{\int_{\Omega_+}p_1^2}\sum_{T\in\mathcal T_+}\int_T
 \begin{pmatrix}
 \nabla p_1^TA\nabla p_a-rbp_1p_a\\
 \nabla p_1^TA^2\nabla p_1\\
 b^2p_1^2\\bp_1p_a
 \end{pmatrix}
 <\begin{pmatrix}-1941428\\43700853\\8980831\\307201\end{pmatrix},
 \quad a=1,\\[6pt]
 \begin{pmatrix}-5269574\\31266387\\1733164\\-685749\end{pmatrix}
 <\frac{10^6}{\int_{\Omega_+}p_1^2}\sum_{T\in\mathcal T_+}\int_T
 \begin{pmatrix}
 \nabla p_1^TA\nabla p_a-rbp_1p_a\\
 \nabla p_1^TA^2\nabla p_1\\
 b^2p_1^2\\bp_1p_a
 \end{pmatrix}
 <\begin{pmatrix}-5269573\\31266388\\1733165\\-685748\end{pmatrix},
 \quad a=2.
 \end{gathered}
\end{equation}
The vector inequalities are componentwise. On each tetrahedron,
\[
 \left(\max_{1\leq i\leq3}\sum_{j=1}^3|A_{ij}|,|b|\right)
 \leq\begin{cases}(10,8),&a=1,\\(6,8/3),&a=2.\end{cases}
\]
The shear test is invariant under $1\leftrightarrow2$. For \eqref{lo:transfer}, the integral bounds give
\[
 \begin{pmatrix}\lVert A\nabla P_1\rVert_2\\\lVert bP_1\rVert_2\end{pmatrix}
 <\begin{cases}
 \begin{pmatrix}331/50\\3\end{pmatrix},&a=1,\\[3pt]
 \begin{pmatrix}6\\7/5\end{pmatrix},&a=2.
 \end{cases}
\]
For $a=2$ the same bounds hold with $P_2$. The first components yield
\[
 (L_r(\diag(1,-1,0)))_{11}<-97/50,\qquad
 (L_r(e_1e_2^T+e_2e_1^T))_{12}<-263/50.
\]
For the diagonal and shear pairings, respectively,
\begin{align*}
 \left|\int_Wbu_1^2\right|
 &\leq\left|\int_WbP_1^2\right|
       +\frac2{125}\lVert bP_1\rVert_2+\frac8{15625}\\
 &<\frac{307201}{10^6}+\frac6{125}+\frac8{15625}\qquad(a=1),\\
 \left|\int_Wbu_1u_2\right|
 &\leq\lVert b\rVert_\infty\lVert u_1\rVert_2\lVert u_2\rVert_2
 \leq8/3\qquad(a=2).
\end{align*}
Equation~\eqref{lo:transfer} gives
\[
 (L_\lambda(\diag(1,-1,0)))_{11}<-7/5,\qquad
 (L_\lambda(e_1e_2^T+e_2e_1^T))_{12}<-24/5.
\]

Coordinate reflection gives
\begin{equation}\label{lo:coefficientidentity}
 \begin{aligned}
 (B_m(K))_{ij}={}&\frac1{8\int_{\Omega_+}p_1^2}
 \sum_{s\in\{-1,1\}^3}\sum_{T\in\mathcal T_+}s_i s_j\\
 &\quad\times\int_T\Bigl(\nabla p_i^TC_m((s_as_bK_{ab})_{a,b=1}^3)\nabla p_j\\
 &\hspace{63pt}{}-r[z^m]\det\!\left(I+zE_T((s_as_bK_{ab})_{a,b=1}^3)\right)p_ip_j\Bigr),
 \end{aligned}
\end{equation}
where the coefficient fields use $E_T((s_as_bK_{ab})_{a,b=1}^3)$ from \eqref{lo:transportmatrix}.
For a distinguished index $d_0$ and $a<b$ with $\{a,b,d_0\}=\{1,2,3\}$, use
\begin{align}\label{lo:Xcoordinates}
 X&=(K_{d_0d_0},K_{aa}-K_{bb},K_{23},K_{13},K_{12}),\\
 K_{aa}&=(-X_1+X_2)/2,\qquad K_{bb}=(-X_1-X_2)/2,\notag\\
 36\lVert p-p_*\rVert_2^2&=12X_1^2+4X_2^2+8\sum_{i=3}^5X_i^2.\notag
\end{align}
Reflections $K_{ij}\mapsto s_is_jK_{ij}$ and permutations imply, for $\{i,j,k\}=\{1,2,3\}$,
\begin{align}\label{lo:quadraticforms}
 (B_2)_{ii}&=\alpha K_{ii}^2+\beta(K_{jj}-K_{kk})^2
       +\gamma_2(K_{ij}^2+K_{ik}^2)+\delta_2K_{jk}^2,\\
 (B_2)_{ij}&=\eta K_{kk}K_{ij}+\zeta K_{ik}K_{jk}\qquad(i\ne j).\notag
\end{align}
In \eqref{lo:coefficientidentity}, extract $(X_1^2,X_2^2,X_4^2,X_3^2)$ from entry $11$ with $d_0=1$, and $(X_1X_5,X_3X_4)$ from entry $12$ with $d_0=3$. The resulting enclosures are
\begin{equation}\label{lo:quadraticintegers}
 \begin{aligned}
 (4853,2868,5418,2931,-9533,2691)
 &<1000(\alpha,\beta,\gamma_2,\delta_2,\eta,\zeta)\\
 &<(4854,2869,5419,2932,-9532,2692).
 \end{aligned}
\end{equation}
For $|\xi|=1$ and $\{i,j,k\}=\{1,2,3\}$,
\begin{align*}
 &\sum_{i=1}^3\xi_i^2\bigl(\alpha K_{ii}^2+\beta(K_{jj}-K_{kk})^2\bigr)\\
 &\quad=2\beta\sum_{i=1}^3K_{ii}^2
       +(\alpha-3\beta)\sum_{i=1}^3K_{ii}^2\xi_i^2
 \leq2\beta\sum_{i=1}^3K_{ii}^2,\\
 &(K_{23}\xi_1+K_{13}\xi_2+K_{12}\xi_3)^2\\
 &\quad\leq3(K_{23}^2\xi_1^2+K_{13}^2\xi_2^2+K_{12}^2\xi_3^2),\\
 &(K_{23}\xi_1+K_{13}\xi_2+K_{12}\xi_3)^2
 \leq\sum_{1\le i<j\le3}K_{ij}^2,\\
 &\delta_2-\gamma_2-\zeta<0,\qquad\delta_2-\gamma_2+2\zeta>0.
\end{align*}
Hence
\begin{align*}
 &\gamma_2\sum_{1\le i<j\le3}K_{ij}^2
       +\zeta(K_{23}\xi_1+K_{13}\xi_2+K_{12}\xi_3)^2\\
 &\qquad+(\delta_2-\gamma_2-\zeta)
       (K_{23}^2\xi_1^2+K_{13}^2\xi_2^2+K_{12}^2\xi_3^2)\\
 &\quad\leq\gamma_2\sum_{1\le i<j\le3}K_{ij}^2
       +\frac{\delta_2-\gamma_2+2\zeta}{3}
        (K_{23}\xi_1+K_{13}\xi_2+K_{12}\xi_3)^2\\
 &\quad\leq\frac{2\gamma_2+\delta_2+2\zeta}{3}
             \sum_{1\le i<j\le3}K_{ij}^2.
\end{align*}
If $\sum_{a=1}^3K_{aa}^2>0$, then
\[
 K_{ii}^2=(K_{jj}+K_{kk})^2\leq2(K_{jj}^2+K_{kk}^2),\qquad
 \frac{K_{ii}^2}{\sum_{a=1}^3K_{aa}^2}\leq\frac23.
\]
Ordering the coordinate squares of $\xi$ as $r_1\geq r_2\geq r_3$ gives
\[
 \frac{K_{11}^2\xi_2^2\xi_3^2+K_{22}^2\xi_1^2\xi_3^2+K_{33}^2\xi_1^2\xi_2^2}
      {\sum_{a=1}^3K_{aa}^2}
 \leq\frac{2r_1r_2+r_1r_3}{3}
 \leq\frac23r_1(1-r_1)\leq\frac16.
\]
The mixed term satisfies
\begin{align*}
 &2|\eta|\left|K_{11}K_{23}\xi_2\xi_3
       +K_{22}K_{13}\xi_1\xi_3+K_{33}K_{12}\xi_1\xi_2\right|\\
 &\quad\leq2|\eta|\left(\sum_{1\le i<j\le3}K_{ij}^2\right)^{1/2}
       \left(K_{11}^2\xi_2^2\xi_3^2+K_{22}^2\xi_1^2\xi_3^2
              +K_{33}^2\xi_1^2\xi_2^2\right)^{1/2}\\
 &\quad\leq|\eta|\left(\frac23
       \left(\sum_{i=1}^3K_{ii}^2\right)\left(\sum_{1\le i<j\le3}K_{ij}^2\right)\right)^{1/2}.
\end{align*}
Consequently,
\begin{align*}
 \xi^TB_2\xi
 &\leq2\beta\sum_{i=1}^3K_{ii}^2
       +\frac{2\gamma_2+\delta_2+2\zeta}{3}\sum_{1\le i<j\le3}K_{ij}^2\\
 &\qquad+|\eta|\left(\frac23
       \left(\sum_{i=1}^3K_{ii}^2\right)\left(\sum_{1\le i<j\le3}K_{ij}^2\right)\right)^{1/2}.
\end{align*}
By \eqref{lo:quadraticintegers},
\[
 10-2\beta>0,\qquad
 \det\left(10I-
 \begin{pmatrix}
 2\beta&|\eta|/\sqrt6\\
 |\eta|/\sqrt6&(2\gamma_2+\delta_2+2\zeta)/3
 \end{pmatrix}\right)>0.
\]
Thus
\[
 \xi^TB_2\xi<10\left(\sum_{i=1}^3K_{ii}^2
                     +\sum_{1\le i<j\le3}K_{ij}^2\right)
 =45\lVert p-p_*\rVert_2^2.
\]

The same symmetries give
\begin{align}\label{lo:cubicforms}
 (B_3)_{11}={}&a_1K_{11}^3+a_2K_{11}(K_{22}-K_{33})^2\notag\\
 &+a_3K_{11}(K_{12}^2+K_{13}^2)+a_4K_{11}K_{23}^2\notag\\
 &+a_5(K_{22}-K_{33})(K_{12}^2-K_{13}^2)
       +a_6K_{12}K_{13}K_{23},\\
 (B_3)_{12}={}&b_1K_{33}^2K_{12}+b_2(K_{11}-K_{22})^2K_{12}\notag\\
 &+b_3K_{33}K_{13}K_{23}+b_4K_{12}^3
       +b_5K_{12}(K_{13}^2+K_{23}^2).\notag
\end{align}
Use $d_0=1$ in entry $11$, $d_0=3$ in entry $12$, and extract the first displayed monomial of each factor ($X_2X_5^2$ for $a_5$). Summation gives
\begin{align}\label{lo:cubicintegers}
 1000(|a_1|,\ldots,|a_6|)&<(4430,6415,15209,7801,15233,3184),\\
 1000(|b_1|,\ldots,|b_5|)&<(26526,223,4934,8108,951).\notag
\end{align}
For each cubic entry, the multinomial formula and Cauchy--Schwarz give
\begin{align*}
 46656\lVert p-p_*\rVert_2^6&=\sum_{\substack{\nu\in\Nzero^5\\|\nu|=3}}
 \binom3\nu12^{\nu_1}4^{\nu_2}8^{\nu_3+\nu_4+\nu_5}X^{2\nu},\\
 \left|\sum_{\substack{\nu\in\Nzero^5\\|\nu|=3}}a_\nu X^\nu\right|^2
 &\leq46656\lVert p-p_*\rVert_2^6\sum_{\substack{\nu\in\Nzero^5\\|\nu|=3}}
 \frac{a_\nu^2}{\binom3\nu12^{\nu_1}4^{\nu_2}8^{\nu_3+\nu_4+\nu_5}}.
\end{align*}
The sums for entries $11$ and $12$ satisfy, respectively,
\[
 \frac1{512}\left(\frac8{27}a_1^2+\frac89a_2^2+\frac49a_3^2
 +\frac29a_4^2+\frac43a_5^2+\frac16a_6^2\right)
 <\frac{459}{500}.
\]
\[
 \frac1{512}\left(\frac4{27}b_1^2+\frac43b_2^2+\frac19b_3^2+b_4^2+\frac23b_5^2\right)
 <\frac{339}{1000}.
\]
Permutation symmetry gives
\[
 \begin{aligned}
 \lVert B_3\rVert_F^2
 &=\sum_{i=1}^3(B_3)_{ii}^2+2\sum_{1\leq i<j\leq3}(B_3)_{ij}^2\\
 &<\left(\frac{2376}{5}\right)^2\lVert p-p_*\rVert_2^6.
 \end{aligned}
\]
\end{proof}

\begin{proof}[Proof of Proposition~\ref{prop:local}]
Since $\lambda_{\min}(Q)\leq-\lVert Q\rVert_F/\sqrt6$ for $Q=Q^T$ with $\tr Q=0$, equations~\eqref{lo:strain} and \eqref{lo:splitting} give
\[
 \begin{aligned}
 \lambda_{\min}(L_\theta(K/2))
 &\leq-\frac{\lVert L_\theta(K/2)\rVert_F}{\sqrt6}\\
 &\leq-\frac{\sqrt3}{4}\lVert p-p_*\rVert_2
 \min\left\{\begin{array}{l}
 \left|(L_\theta(\diag(1,-1,0)))_{11}\right|,\\
 \sqrt2\left|(L_\theta(e_1e_2^T+e_2e_1^T))_{12}\right|
 \end{array}\right\}.
 \end{aligned}
\]
Numerical Lemma~\ref{num:local} gives
\begin{equation}\label{lo:linearnegative}
 \lambda_{\min}(L_\lambda(K/2))<-\frac{291}{500}\lVert p-p_*\rVert_2,
 \qquad\lambda_{\min}(L_r(K/2))<-\frac{21}{25}\lVert p-p_*\rVert_2.
\end{equation}
Using \eqref{lo:transfer}, $\lVert\nabla P_\xi\rVert_2<21/10$, and
\eqref{lo:homogeneousbounds}, the bounds \eqref{lo:highercoefficients} give
\[
 \begin{aligned}
 \mathcal B_\lambda\!\left(u_\xi,u_\xi;C_2,
          \frac{(\tr E)^2-\tr(E^2)}2\right)&<\frac{153}{2}\lVert p-p_*\rVert_2^2,\\
 \mathcal B_\lambda(u_\xi,u_\xi;C_3,\det E)&<756\lVert p-p_*\rVert_2^3.
 \end{aligned}
\]

\begin{align*}
 C-\sum_{m=0}^3C_m&=J\sum_{m=4}^{\infty}M_m+(J-1)M_3\\
 &\quad+\left(\frac{(\tr E)^2-\tr(E^2)}2+\det E\right)M_2+(\det E)M_1,\\
 \left\lVert\sum_{m=4}^{\infty}M_m\right\rVert_{\rm op}
 &\leq\sum_{m=4}^{\infty}(1+3m/2)(6\lVert p-p_*\rVert_2)^m\\
 &=\frac{1296\lVert p-p_*\rVert_2^4(7-33\lVert p-p_*\rVert_2)}
        {(1-6\lVert p-p_*\rVert_2)^2}\\
 &\leq\frac{9072\lVert p-p_*\rVert_2^4}{(1-6\lVert p-p_*\rVert_2)^2}.
\end{align*}
Since
\[
 |J-1|\leq\frac{26}{5}\lVert p-p_*\rVert_2
 +\frac{63}{2}\lVert p-p_*\rVert_2^2
 +\frac{216}{5}\lVert p-p_*\rVert_2^3,
\]
\begin{align*}
 &\frac{\left\lVert C-\sum_{m=0}^3C_m\right\rVert_{\rm op}}
       {1296\lVert p-p_*\rVert_2^4}\\
 &\quad\leq\frac{7(1+|J-1|)}{(1-6\lVert p-p_*\rVert_2)^2}
 +\frac{11|J-1|}{12\lVert p-p_*\rVert_2}
 +4\left(\frac78+\frac65\lVert p-p_*\rVert_2\right)+\frac12\\
 &\quad<\begin{cases}
 35/2,&\lVert p-p_*\rVert_2\leq1/120,\\
 37/2,&\lVert p-p_*\rVert_2\leq1/75.
 \end{cases}
\end{align*}
For each unit mode $f$,
\begin{align*}
 &\left|\int_W\nabla f^T\left(C-\sum_{m=0}^3C_m\right)\nabla f\right|\\
 &\quad\leq\frac{430113}{100000}
       \left\lVert C-\sum_{m=0}^3C_m\right\rVert_{\rm op}\\
 &\quad<\begin{cases}
 98496\lVert p-p_*\rVert_2^4,&\lVert p-p_*\rVert_2\leq1/120,\\
 103680\lVert p-p_*\rVert_2^4,&\lVert p-p_*\rVert_2\leq1/75.
 \end{cases}
\end{align*}

For $0<\lVert p-p_*\rVert_2\leq27/4000$, choose $u_\xi$ to minimize
the linear term. Then
\begin{equation}\label{lo:I}
 \begin{aligned}
 \mathcal B_\lambda(u_\xi,u_\xi;C,J)
 &<-\frac{291}{500}\lVert p-p_*\rVert_2
   +\frac{153}{2}\lVert p-p_*\rVert_2^2\\
 &\quad+756\lVert p-p_*\rVert_2^3
       +98496\lVert p-p_*\rVert_2^4<0.
 \end{aligned}
\end{equation}
The quotient by $\lVert p-p_*\rVert_2$ is increasing and negative at
$\lVert p-p_*\rVert_2=27/4000$.
For $27/4000\leq\lVert p-p_*\rVert_2\leq73/6000$, choose $P_\xi$ to
minimize its linear term. Since
\[
 \begin{gathered}
 \mathcal B_{\frac{4300827}{10^6}}(P_\xi,P_\xi;C,J)
 =\mathcal B_r(P_\xi,P_\xi;C,J)
 +(r-\tfrac{4300827}{10^6})\int_WJP_\xi^2,\\
 \int_WJP_\xi^2\leq1+\frac{57}{10}\lVert p-p_*\rVert_2,
 \end{gathered}
\]
we obtain
\begin{equation}\label{lo:A}
 \begin{aligned}
 \mathcal B_{\frac{4300827}{10^6}}(P_\xi,P_\xi;C,J)
 &<-\frac{21}{25}\lVert p-p_*\rVert_2+45\lVert p-p_*\rVert_2^2\\
 &\quad+\frac{2376}{5}\lVert p-p_*\rVert_2^3
       +103680\lVert p-p_*\rVert_2^4\\
 &\quad+\frac{303}{10^6}
       \left(1+\frac{57}{10}\lVert p-p_*\rVert_2\right)<0.
 \end{aligned}
\end{equation}
The last polynomial is convex and negative at both endpoints.
Thus \eqref{lo:Rayleigh} gives
\begin{equation}\label{lo:finalcomparison}
 \begin{gathered}
 2^{-4/3}F(p)<
 \begin{cases}
 \lambda,&0<\lVert p-p_*\rVert_2\leq27/4000,\\
 \frac{4300827}{10^6},&27/4000<\lVert p-p_*\rVert_2\leq73/6000,
 \end{cases}\\
 \frac{4300827}{10^6}<\lambda=2^{-4/3}F(p_*).
 \end{gathered}
\end{equation}
\end{proof}

\section{The global comparison}\label{sec:global}
Geometric exclusions and polynomial Rayleigh tests give the remaining comparison:
\begin{equation}\label{gl:target}
 D(p)>0,\quad\norm{p-p_*}_2>\frac{73}{6000}
 \quad\Longrightarrow\quad F(p)<\frac{10837}{1000}=T<F(p_*).
\end{equation}
Together with Proposition~\ref{prop:local}, this proves Theorem~\ref{thm:main}.
Relabel an edge of largest weight as $01$. Its stabilizer acts transitively
on $02,03,12,13$; permutation and rescaling give
\[
 F(sp)=F(p)\quad(s>0),\qquad
 q_e=\frac{p_e}{p_{01}},\qquad q_{02}=\max(q_{02},q_{03},q_{12},q_{13}).
\]
It therefore suffices to consider
\begin{equation}\label{gl:chamber}
 \begin{gathered}
 \mathcal C=\{q:q_{01}=1,\ 0\le q_e\le1\ (e\in E),\
 q_{02}\ge\max(q_{03},q_{12},q_{13})\}.
 \end{gathered}
\end{equation}
We evaluate $F,D,C_e,W_S$ at these unnormalized conorms.
For $w=w_1\cdots w_h\in\{0,1\}^h$, bisect $02,03,12,13,23$ cyclically,
choosing the lower half for $0$ and the upper half for $1$:
\[
 \begin{gathered}
 B_w=\{1\}\times\prod_{i=1}^5 2^{-k_i}[j_i,j_i+1],\\
 k_i=\#\big(\{i,i+5,\ldots\}\cap\{1,\ldots,h\}\big),\qquad
 j_i=\sum_{r=0}^{k_i-1}2^{k_i-1-r}w_{i+5r}.
 \end{gathered}
\]
Here $k_i=j_i=0$ for an empty subsequence; $q^-,q^+$ denote the box endpoints. We have
\begin{equation}\label{gl:emptytest}
 q_{02}^+<\max(q_{03}^-,q_{12}^-,q_{13}^-)
 \quad\Longrightarrow\quad B_w\cap\mathcal C=\varnothing.
\end{equation}
\begin{equation}\label{gl:localtest}
 \begin{gathered}
 \max_{v\in\operatorname{Vert}(B_w)}
 \left(\frac{\sum_{0\le i<j\le3}v_{ij}^2}{(\sum_{0\le i<j\le3}v_{ij})^2}-\frac16\right)
 \le\left(\frac{73}{6000}\right)^2\\
 \Longrightarrow\quad
 \left\lVert\frac{q}{\sum_{0\le i<j\le3}q_{ij}}-p_*\right\rVert_2
 \le\frac{73}{6000}\qquad(q\in B_w).
 \end{gathered}
\end{equation}
Indeed, for $q=\sum_{v\in\operatorname{Vert}(B_w)}a_vv$, $a_v\ge0$, $\sum_{v\in\operatorname{Vert}(B_w)}a_v=1$,
\[
 \frac{q}{\sum_{0\le i<j\le3}q_{ij}}=\sum_{v\in\operatorname{Vert}(B_w)}
 \frac{a_v\sum_{0\le i<j\le3}v_{ij}}{\sum_{0\le i<j\le3}q_{ij}}\frac{v}{\sum_{0\le i<j\le3}v_{ij}}.
\]
Convexity of the squared norm therefore gives
\[
 \left\lVert\frac{q}{\sum_{0\le i<j\le3}q_{ij}}-p_*\right\rVert_2^2
 =\frac{\sum_{0\le i<j\le3}q_{ij}^2}
        {(\sum_{0\le i<j\le3}q_{ij})^2}-\frac16
 \le\max_{v\in\operatorname{Vert}(B_w)}\left\lVert\frac{v}{\sum_{0\le i<j\le3}v_{ij}}-p_*\right\rVert_2^2.
\]

The sixteen spanning trees of $K_4$, in lexicographic order, are
\begin{equation}\label{gl:treelist}
\begin{array}{llll}
 J_1=\{01,02,03\},&J_2=\{01,02,13\},&J_3=\{01,02,23\},&J_4=\{01,03,12\},\\
 J_5=\{01,03,23\},&J_6=\{01,12,13\},&J_7=\{01,12,23\},&J_8=\{01,13,23\},\\
 J_9=\{02,03,12\},&J_{10}=\{02,03,13\},&J_{11}=\{02,12,13\},&J_{12}=\{02,12,23\},\\
 J_{13}=\{02,13,23\},&J_{14}=\{03,12,13\},&J_{15}=\{03,12,23\},&J_{16}=\{03,13,23\}.
\end{array}
\end{equation}
Cauchy--Binet and the matrix--tree identity~\cite{Chaiken} give
\begin{equation}\label{gl:trees}
 \begin{aligned}
 D(q)&=\sum_{\substack{J\subset E\\\#J=3}}
       \bigl(\det[d_e]_{e\in J}\bigr)^2\prod_{e\in J}q_e
       =\sum_{r=1}^{16}\prod_{e\in J_r}q_e,\\
 C_{ij}(q)&=\partial_{q_{ij}}D(q)
 =d_{ij}^T\adj(G(q))d_{ij}\quad(0\le i<j\le3).
 \end{aligned}
\end{equation}
By \eqref{eq:dualbound}, \eqref{eq:primalbound}, and \eqref{eq:cornerconstants},
\begin{equation}\label{gl:cornertest}
 \begin{gathered}
 \frac{C_e(q)^3}{D(q)^2}<\frac{132379}{100000}
 \quad\text{or}\quad\frac{W_S(q)^3}{D(q)}<\frac{28}{25},\\
 D(q)>0\quad\Longrightarrow\quad F(q)<T.
 \end{gathered}
\end{equation}
For $a,b,c,d,z\ge0$ and $c+dz>0$,
\[
 \frac{\dd^2}{\dd z^2}\frac{(a+bz)^3}{(c+dz)^m}\ge0,
 \qquad m\in\{1,2\}.
\]
Their maxima are therefore attained at corners:
\[
 \max_{q\in B_w}\frac{C_e(q)^3}{D(q)^2}
 =\max_{v\in\operatorname{Vert}(B_w)}\frac{C_e(v)^3}{D(v)^2},\qquad
 \max_{q\in B_w}\frac{W_S(q)^3}{D(q)}
 =\max_{v\in\operatorname{Vert}(B_w)}\frac{W_S(v)^3}{D(v)}.
\]
For cut tests require $\min_{v\in\operatorname{Vert}(B_w)}D(v)>0$; for edge tests require
$q_e^->0$ and use the continuous extension by zero at $D=0$:
\[
 D(q)=D(q)|_{q_e=0}+q_eC_e(q)\ge q_eC_e(q),\qquad
 0\le\frac{C_e^3}{D^2}\le\frac{D}{(q_e^-)^3}\longrightarrow0.
\]

For an odd trial in \eqref{eq:functional}, let $M,K$ be its mass and
energy after canceling their common positive Jacobian. For the box center $c$,
\[
 \begin{gathered}
 F(q)\le D(q)^{1/3}\frac{K(q)}{M(q)}\quad(M>0),\\
 \sum_{0\le i<j\le3}C_{ij}(c)q_{ij}
       =\tr\{\adj(G(c))G(q)\},\qquad\sum_{0\le i<j\le3}C_{ij}(c)c_{ij}=3D(c).
 \end{gathered}
\]
The arithmetic--geometric mean inequality now gives
\[
 \left(\frac{D(q)}{D(c)}\right)^{1/3}
 \le\frac13\tr\{G(c)^{-1/2}G(q)G(c)^{-1/2}\}
 =\frac{\sum_{0\le i<j\le3}C_{ij}(c)q_{ij}}{3D(c)},
\]
so either test
\begin{align}
 0<\tau^3<27T^3D(c)^2,\qquad
 &\tau M-\left(\sum_{0\le i<j\le3}C_{ij}(c)q_{ij}\right)K>0,\label{gl:linearRayleigh}\\
 M>0,\qquad &T^3M^3-DK^3>0\label{gl:cubedRayleigh}
\end{align}
gives $F(q)<T$. For \eqref{gl:linearRayleigh},
\[
 K\ge0\quad\Longrightarrow\quad M>0,\qquad
 F(q)\le\frac{\left(\sum_{0\le i<j\le3}C_{ij}(c)q_{ij}\right)K}{3D(c)^{2/3}M}
 <\frac{\tau}{3D(c)^{2/3}}<T.
\]
Order the edges as $(01,02,03,12,13,23)$ and
the trees as in \eqref{gl:treelist}, with heights
\[
 (h_{01},h_{02},h_{03},h_{12},h_{13},h_{23})=(1,8,64,512,4096,32768).
\]
For each $J\in\{J_1,\ldots,J_{16}\}$, solve
\[
 d_e\cdot\ell_J=h_e\qquad(e\in J).
\]
The signed sums of distinct tree heights satisfy
$h_e-d_e\cdot\ell_J\ne0$ for $e\notin J$. For $q>0$, the parallelotopes
\begin{equation}\label{gl:tiles}
 \frac12\sum_{e\in E\setminus J}\sgn(h_e-d_e\cdot\ell_J)q_ed_e
 +\sum_{e\in J}[-q_ed_e/2,q_ed_e/2]
\end{equation}
tile $Y(q)$ with disjoint interiors.

Indeed, under $\sum_{0\le i<j\le3}z_{ij}d_{ij}=y$, $\abs{z_e}\le q_e/2$,
\[
 \begin{split}
 \sum_{0\le i<j\le3}h_{ij}z_{ij}
 &=\ell_J\cdot y+\sum_{e\in E\setminus J}(h_e-d_e\cdot\ell_J)z_e\\
 &\le\ell_J\cdot y+\frac12\sum_{e\in E\setminus J}
                           \abs{h_e-d_e\cdot\ell_J}q_e,
 \end{split}
\]
with equality precisely for $z_e=\sgn(h_e-d_e\cdot\ell_J)q_e/2$ off $J$.
At a vertex maximizer, the free columns are independent and, for almost
every $y$, form a tree $J$. A wrong off-tree sign permits, for small $t>0$,
\[
 \delta z_e=t\sgn(h_e-d_e\cdot\ell_J),\qquad
 \delta z_J=-[d_a]_{a\in J}^{-1}d_e\,\delta z_e,\qquad
 \delta z_a=0\quad(a\in E\setminus(J\cup\{e\})),
\]
so that
\[
 \sum_{0\le i<j\le3}d_{ij}\delta z_{ij}=0,\qquad
 \sum_{0\le i<j\le3}h_{ij}\delta z_{ij}=t\abs{h_e-d_e\cdot\ell_J}>0,
\]
contradicting maximality. The off-tree signs and the tree coordinates are
unique; hence the interiors are disjoint, and closure gives \eqref{gl:tiles}.

In the fixed coordinates of Section~\ref{sec:local},
\begin{equation}\label{gl:fixedframe}
 \begin{gathered}
 \Gamma(q)=\frac14
 \begin{pmatrix}
 q_{02}+q_{03}+q_{12}+q_{13}&q_{03}-q_{12}&q_{02}-q_{13}\\
 q_{03}-q_{12}&q_{01}+q_{03}+q_{12}+q_{23}&q_{01}-q_{23}\\
 q_{02}-q_{13}&q_{01}-q_{23}&q_{01}+q_{02}+q_{13}+q_{23}
 \end{pmatrix},\\
 x_J(z)=\frac14\sum_{0\le i<j\le3}q_{ij}(b_i-b_j)z_{ij},
 \qquad z_e=\sgn(h_e-d_e\cdot\ell_J)/2\quad(e\notin J).
 \end{gathered}
\end{equation}
Since $\abs{\det[d_e]_{e\in J_r}}=1$, $1\le r\le16$, canceling the
common frame Jacobian gives
\begin{equation}\label{gl:integral}
 \begin{split}
 M&=\sum_{r=1}^{16}\left(\prod_{e\in J_r}q_e\right)
       \int_{[-1/2,1/2]^3}f(x_{J_r}(z))^2\dd z_{J_r},\\
 K&=\sum_{r=1}^{16}\left(\prod_{e\in J_r}q_e\right)
       \int_{[-1/2,1/2]^3}
       (\nabla f^T\Gamma\nabla f)(x_{J_r}(z))\dd z_{J_r}.
 \end{split}
\end{equation}
These polynomial expressions extend to $q_e=0$, using
\[
 \int_{-1/2}^{1/2}z^{2j}\dd z=\frac{2^{-2j}}{2j+1},\qquad
 \int_{-1/2}^{1/2}z^{2j+1}\dd z=0\quad(j\in\Nzero).
\]

The one-dimensional rules are
\[
 \begin{gathered}
 Q_{-1}f=0,\qquad Q_0f=f(0),\\
 Q_jf=\sum_{k=-j}^j f\!\left(\frac{k}{2j}\right)
 \int_{-1/2}^{1/2}
 \prod_{\substack{\ell\in\Z\\-j\le\ell\le j,\ \ell\ne k}}
       \frac{2jx-\ell}{k-\ell}\,\dd x\qquad(j\ge1).
 \end{gathered}
\]
For each $s\in\{3,5\}$, add coincident rational weights before rounding and define
$w_y$, zero at all other points, by
\begin{equation}\label{gl:cubature}
 \sum_{\substack{y\in[-1/2,1/2]^3\\w_y\ne0}}w_y f(y)
 =\sum_{i=0}^s\sum_{j=0}^{s-i}\sum_{k=0}^{s-i-j}
 [(Q_i-Q_{i-1})\otimes(Q_j-Q_{j-1})\otimes(Q_k-Q_{k-1})]f.
\end{equation}
Telescoping gives, for $a,b,c\in\Nzero$ with $a+b+c\le s$,
\[
 \sum_{\substack{y\in[-1/2,1/2]^3\\w_y\ne0}}
 w_y y_1^{2a}y_2^{2b}y_3^{2c}
 =\frac{2^{-2(a+b+c)}}{(2a+1)(2b+1)(2c+1)},
\]
\[
 \#\{y\in[-1/2,1/2]^3:w_y\ne0\}
 =\begin{cases}69,&s=3,\\357,&s=5.\end{cases}
\]
Monomials with an odd exponent vanish by symmetry.

The cubic monomial order is
\[
 x_1,x_2,x_3,x_1^3,x_2^3,x_3^3,
 x_1^2x_2,x_1^2x_3,x_1x_2^2,x_2^2x_3,x_1x_3^2,x_2x_3^2,x_1x_2x_3.
\]
For quintics, order $x_1^ax_2^bx_3^c$, $a+b+c\in\{1,3,5\}$, by
increasing $(a+b+c,a,b)$.

Apply coefficientwise truncation
\[
 \mathcal T_h(a)=h\,\sgn(a)\left\lfloor\frac{\abs a}{h}\right\rfloor.
\]
For FCC, the box center $c$ determines the coordinates
\begin{equation}\label{gl:eta}
 (z_1,z_2,z_3)=
 2^{-24}\left\lfloor
 \frac{12\cdot2^{24}}{3(c_{01}+c_{02}+c_{13}+c_{23})+4(c_{03}+c_{12})}\right\rfloor
 (x_3,x_1+x_2,x_1-x_2).
\end{equation}
Define
\begin{equation}\label{gl:fccpoly}
 \begin{split}
 10^6g(z)=z_1\big[&10^6-860262(z_2^2+z_3^2)
 -109471(z_2^4+z_3^4)+3253797z_2^2z_3^2\\
 &-z_1^2\{460101+291610(z_2^2+z_3^2)\}+99272z_1^4\big].
 \end{split}
\end{equation}
Interchanging $z_1,z_i$ gives $g_i$ and the trial
\[
 \mathcal T_{2^{-40}}\left(\sum_{i=1}^3v_i
                    \mathcal T_{2^{-40}}g_i\right).
\]
Its energy matrix is
\begin{multline}\label{gl:fccmetric}
 2^{-50}\left\lfloor
 \frac{12\cdot2^{24}}{3(c_{01}+c_{02}+c_{13}+c_{23})+4(c_{03}+c_{12})}\right\rfloor^2\\
 {}\times
 \begin{pmatrix}
 q_{01}+q_{02}+q_{13}+q_{23}
 &q_{01}+q_{02}-q_{13}-q_{23}
 &-q_{01}+q_{02}-q_{13}+q_{23}\\
 q_{01}+q_{02}-q_{13}-q_{23}
 &q_{01}+q_{02}+4q_{03}+q_{13}+q_{23}
 &-q_{01}+q_{02}+q_{13}-q_{23}\\
 -q_{01}+q_{02}-q_{13}+q_{23}
 &-q_{01}+q_{02}+q_{13}-q_{23}
 &q_{01}+q_{02}+4q_{12}+q_{13}+q_{23}
 \end{pmatrix}.
\end{multline}

For transported trials, use the common BCC fan \eqref{sp:tetrahedra}, with its
fixed tetrahedron order and coordinate reflections.
Let $F_A(q)$ send each reference vertex and face average to its
corresponding point in the zonotope with weights $q_e$. The vertex of
strict order $\sigma$ is
\begin{equation}\label{gl:vertex}
 V_q(\sigma)=\frac18\sum_{0\le i<j\le3}q_{ij}
           \sgn(\sigma_i-\sigma_j)(b_i-b_j).
\end{equation}
Positive segment lengths preserve the normal fan, and face averages
remain in relative interiors. Hence
\[
 q_e>0\ (e\in E)\quad\Longrightarrow\quad
 \text{the face incidences persist},\quad
 \det F_A(q)>0.
\]
The coned triangulations are therefore nondegenerate.

Transporting the normalized modes of Numerical Lemma~\ref{num:local}
preserves oddness and gives
\begin{equation}\label{gl:transportforms}
 \begin{aligned}
 M&=\sum_{T\in\mathcal T_+}\sum_{s\in\{-1,1\}^3}
       \det(F_{\diag(s_1,s_2,s_3)T})
       \int_{\diag(s_1,s_2,s_3)T}\left(\sum_{i=1}^3v_iP_i\right)^2,\\
 K&=\sum_{T\in\mathcal T_+}\sum_{s\in\{-1,1\}^3}
       \frac1{\det F_{\diag(s_1,s_2,s_3)T}}
       \int_{\diag(s_1,s_2,s_3)T}
       \left(\sum_{i=1}^3v_i\nabla P_i\right)^T\\
  &\qquad\qquad\times\adj(F_{\diag(s_1,s_2,s_3)T})
       \Gamma\adj(F_{\diag(s_1,s_2,s_3)T})^T\\
  &\qquad\qquad\times
       \left(\sum_{j=1}^3v_j\nabla P_j\right).
 \end{aligned}
\end{equation}
All reference integrals are rational by \eqref{lo:moments}.

The Appendix specifies the center coefficients, enclosure
arithmetic, and test order.
Accept a box at the first successful test, and set
\begin{equation}\label{gl:recursion}
 \mathcal W_0=\{\varnothing\},\qquad
 \mathcal W_{h+1}=\{wb:w\in\mathcal W_h,\ B_w\text{ is not accepted},\ b\in\{0,1\}\}.
\end{equation}

\begin{lemma}\label{gl:covering}
If $\mathcal W_{H+1}=\varnothing$ for some $H$, then every $q\in\mathcal C$
with $D(q)>0$ satisfies
\[
 F(q)<T\qquad\text{or}\qquad\left\lVert\frac{q}{\sum_{0\le i<j\le3}q_{ij}}-p_*\right\rVert_2\le\frac{73}{6000}.
\]
\end{lemma}

\begin{proof}
Since $B_w=B_{w0}\cup B_{w1}$, \eqref{gl:recursion} gives,
by induction, for every $h\ge0$,
\[
 \{1\}\times[0,1]^5=
 \bigcup_{j=0}^{h-1}\ \bigcup_{\substack{w\in\mathcal W_j\\B_w\text{ is accepted}}}B_w
 \ \cup\!\bigcup_{w\in\mathcal W_h}B_w.
\]
If $\mathcal W_{H+1}=\varnothing$, every $q\in\mathcal C$ lies in an
accepted box. Since the chamber test excludes no such $q$,
\[
 \eqref{gl:localtest}\ \Longrightarrow\ \left\lVert\frac{q}{\sum_{0\le i<j\le3}q_{ij}}-p_*\right\rVert_2\le\frac{73}{6000},
 \qquad
 \eqref{gl:cornertest},\ \eqref{gl:linearRayleigh},\ \eqref{gl:cubedRayleigh}
 \ \Longrightarrow\ F(q)<T.
\]
\end{proof}

\begin{numlemma}\label{num:global}
For \eqref{gl:recursion} and the Appendix constructions,
\[
 \#\mathcal W_{37}=4,\qquad\mathcal W_{38}=\varnothing.
\]
\end{numlemma}

The Appendix gives the acceptance counts and strict sign bounds.

\begin{proof}[Proof of Theorem~\ref{thm:main}]
By Lemma~\ref{gl:covering}, Numerical Lemma~\ref{num:global}, and
Proposition~\ref{prop:local}, a reduced representative satisfies
\[
 \begin{gathered}
 q\in\mathcal C,\quad D(q)>0,\quad F(q)\ge F(p_*)>T\\
 \Longrightarrow\quad\left\lVert\frac{q}{\sum_{0\le i<j\le3}q_{ij}}-p_*\right\rVert_2\le\frac{73}{6000}
 \quad\Longrightarrow\quad \frac{q}{\sum_{0\le i<j\le3}q_{ij}}=p_*.
 \end{gathered}
\]
Equal conorms represent $\BCC$; for $R\in O(3)$,
$\mu_1(RV_{\BCC})=\mu_1(V_{\BCC})$.
Finally, the mixed eigenvalue enclosure and
\[
 V_{\BCC}=2^{-2/3}W,\qquad
 \mu_1(V_{\BCC})=2^{4/3}\mu_1(W)
\]
give the stated enclosure.
\end{proof}

\appendix
\section*{Appendix}
\setcounter{equation}{0}
\renewcommand{\theequation}{A.\arabic{equation}}
\renewcommand{\theHequation}{appendix.\arabic{equation}}
We give the coefficient construction for Numerical Lemma~\ref{num:flux}
and the rational enclosures for Numerical Lemma~\ref{num:global}.

\begin{proof}[Proof of Numerical Lemma~\ref{num:flux}]
Use \eqref{sp:tetrahedra} and the Bernstein convention of Section~\ref{sec:spectrum}.
For ordered vertices $w_0,\ldots,w_3$, with lexicographic sorting,
\[
 \begin{aligned}
 k(P,\alpha)&=\operatorname{sort}\{(2w_i,\alpha_i):0\leq i\leq3,\ \alpha_i>0\},\\
 k^*&=\operatorname{sort}\{((u_1,u_3,u_2),a):((u_1,u_2,u_3),a)\in k\}.
 \end{aligned}
\]

\emph{The scalar space.}
Lexicographically index the distinct $\min(k,k^*)$ with $|\alpha|=9$
and support not contained in $\{x_1=0\}$. Then
\[
 (s_i)_{P,\alpha}=
 \begin{cases}
 1,&\min(k(P,\alpha),k(P,\alpha)^*)\text{ has index }i,\\
 0,&\text{otherwise}.
 \end{cases}
 \qquad
 s_i\in C(\overline{\Omega_+}),\quad s_i|_{x_1=0}=0
 \quad(1\le i\le710).
\]
Set
\begin{equation}\label{sp:scalar-matrices}
 M_{ij}=(s_i,s_j),\qquad
 F_{ij}=(\nabla s_i,\nabla s_j)-\frac{107}{25}(s_i,s_j),
\end{equation}
with inner products over $\Omega_+$ and $1\leq i,j\leq710$.

\emph{The flux space.}
For $|\alpha|=10$, group coefficients by $k$, ordered by first occurrence
in $(P,\alpha)$: $P$ follows \eqref{sp:tetrahedra}, and $\alpha$ is lexicographic.
Within each group, order by $(P,\text{component})$. With
$\operatorname{supp}k=\{u/2:(u,a)\in k\}$, impose
\[
 \begin{aligned}
 (\xi_{P,k}-\xi_{P',k})\cdot\bigl((u_b-u_a)\times(u_c-u_a)\bigr)&=0
 &&\bigl(F=P\cap P',\ \operatorname{supp}k\subset F\bigr),\\
 \xi_{P,k}\cdot\bigl((u_b-u_a)\times(u_c-u_a)\bigr)&=0
 &&\bigl(F\subset\partial\Omega_+\text{ Neumann},\
 \operatorname{supp}k\subset F\bigr),
 \end{aligned}
\]
where $u_a,u_b,u_c$ are the lexicographically ordered doubled face vertices.
For each free index $i$, solve the reduced-row-echelon constraints with
free coordinates $(z_i)_j=\delta_{ij}$, then successively
\[
 z_i\leftarrow\left(\mathop{\operatorname{lcm}}_{j\text{ over entries}}
       \operatorname{den}((z_i)_j)\right)z_i,
 \qquad
 z_i\leftarrow\frac{z_i}{\mathop{\gcd}_{j\text{ over entries}}|(z_i)_j|}.
\]
Denominators are positive. Relabel consecutively by $(\text{group},\text{free index})$:
\[
 \{z_i\}_{i=1}^{6138}\subset H(\operatorname{div};\Omega_+),\qquad
 z_i\cdot n=0\quad\hbox{on the Neumann boundary}.
\]

\emph{Rational coefficient generation.}
For scalar $a$, square brackets mean rounding to the nearest finite
IEEE 754 binary64 value:
\[
 [a]\in\operatorname*{argmin}_{\substack{b=m2^e,\ m,e\in\Z\\
                  |m|<2^{53},\ -1074\le e\le971}}|a-b|,
 \qquad [G]_{ij}=[G_{ij}].
\]
At a midpoint choose the even canonical significand; a zero result retains
the input sign. Vector and matrix inputs are rounded entrywise, and every
nested bracket is a separate rounding. Define nearest-integer rounding by
\[
 \operatorname{nint}(a)=
 \begin{cases}
 \lfloor a+\tfrac12\rfloor,&a\notin\Z+\tfrac12,\\
 2\lfloor\tfrac a2+\tfrac14\rfloor,&a\in\Z+\tfrac12.
 \end{cases}
\]
Thus midpoint ties select the even integer. For symmetric $G\in\R^{N\times N}$,
\[
 G^{(1)}_{ij}=[G_{ij}]\qquad(1\le j\le i\le N).
\]
For $k=1,\ldots,N$,
\begin{equation}\label{sp:ldl-recurrence}
 \begin{aligned}
 l_{ik}&=[G^{(k)}_{ik}/G^{(k)}_{kk}] &&(i>k),\\
 G^{(k+1)}_{ij}
 &=[G^{(k)}_{ij}-[G^{(k)}_{ik}l_{jk}]]&&(i\geq j>k).
 \end{aligned}
\end{equation}
For rounded $y$, compute $S_G(y)$ by the ordered substitutions
\[
 \begin{array}{lll}
 k=1,\ldots,N:&y_i\leftarrow[y_i-[l_{ik}y_k]],&i=k+1,\ldots,N,\\[2pt]
 k=1,\ldots,N:&y_k\leftarrow[y_k/G^{(k)}_{kk}],&\\[2pt]
 k=N,\ldots,1:&y_k\leftarrow[y_k-[l_{ik}y_i]],&i=k+1,\ldots,N.
 \end{array}
\]
All pivots below are positive.

For doubled vertices $u_i=2w_i$, initialize $a^{(0)}$ by rounding
\[
 (x_1)_{P,\alpha}=\frac1{18}\sum_{i=0}^3(u_i)_1\alpha_i,
 \qquad \alpha\in\Nzero^4,\quad|\alpha|=9.
\]
For $j=0,\ldots,7$ and $i=1,\ldots,710$, set
\[
 a_i^{(j+1)}=
 \left[\frac{(S_F([M]a^{(j)}))_i}
 {\max_{1\le k\le710}|(S_F([M]a^{(j)}))_k|}\right].
\]
Evaluate the solve and its fixed maximum once per step. In $[M]a^{(j)}$,
round products and partial sums in increasing column order. Set
\begin{equation}\label{sp:trial-generation}
 v=10^{-10}\sum_{i=1}^{710}\operatorname{nint}(10^{10}a_i^{(8)})s_i.
\end{equation}
Set
\begin{equation}\label{sp:flux-generation}
 \begin{aligned}
 H_{ij}&=(z_i,z_j)+20(\operatorname{div}z_i,\operatorname{div}z_j),\\
 q&=10^{-10}\sum_{i=1}^{6138}
 \operatorname{nint}\!\left(10^{10}\left(S_H\!\left(
 \left[((z_j,\nabla v)-20\frac{\norm{\nabla v}_2^2}{\norm v_2^2}(\operatorname{div}z_j,v))_{j=1}^{6138}\right]
 \right)\right)_i\right)z_i.
 \end{aligned}
\end{equation}
Use rational arithmetic for the following integrations.

\emph{Integration.}
On $P$,
\[
 \nabla\zeta_i=[w_1-w_0,w_2-w_0,w_3-w_0]^{-T}e_i\quad(1\le i\le3),
 \qquad\nabla\zeta_0=-\sum_{i=1}^3\nabla\zeta_i.
\]
For degree $m$,
\begin{equation}\label{sp:bernstein-derivative}
 (\partial_jf)_{P,\alpha}
 =m\sum_{i=0}^3(\partial_j\zeta_i)f_{P,\alpha+e_i},
 \qquad \alpha\in\Nzero^4,\quad|\alpha|=m-1,\quad1\le j\le3.
\end{equation}
Equation~\eqref{sp:bernstein-product} evaluates \eqref{sp:scalar-matrices}, \eqref{sp:flux-generation}, and
\[
 (v,v),\quad(\nabla v,\nabla v),\quad(q,q),\quad
 (\nabla v,q),\quad(\operatorname{div}q,\operatorname{div}q).
\]
The boundary conditions give
\[
 (v,\operatorname{div}q)+(\nabla v,q)
 =\int_{\partial\Omega_+}v\,q\cdot n=0.
\]
Evaluation gives
\begin{equation}\label{sp:integer-parts}
 \begin{aligned}
 \left\lfloor10^{10}\min_{1\le k\le710}G^{(k)}_{kk}\right\rfloor
 &=198287\qquad(G=F),\\
 \left\lfloor10^{12}\min_{1\le k\le6138}G^{(k)}_{kk}\right\rfloor
 &=3063\qquad(G=H),\\
 \left\lfloor10^{12}\frac{\norm{\nabla v}_2^2}{\norm v_2^2}\right\rfloor
 &=4300846940707,\\
 \left\lfloor10^{12}\frac{\norm{\nabla v-q}_2^2}{\norm v_2^2}\right\rfloor
 &=9001492,\\
 \left\lfloor\frac{10^{17}}{\norm v_2^2}
   \left\lVert\operatorname{div}q+
    \frac{\norm{\nabla v}_2^2}{\norm v_2^2}v\right\rVert_2^2\right\rfloor
 &=226155.
 \end{aligned}
\end{equation}
\[
 \min_{0\le j\le7}\max_{1\le i\le710}|(S_F([M]a^{(j)}))_i|>43.
\]
The last three integer parts follow from \eqref{sp:bernstein-product} by
cross-multiplication. Thus \eqref{sp:trial-bounds} holds.

\end{proof}

\begin{proof}[Proof of Numerical Lemma~\ref{num:global}]
In the five free coordinates,
\[
 q_e=c_e+r_eu_e,\qquad |u_e|\le1.
\]
The enclosure $(P,r)$ satisfies
\[
 \abs{h(u)-P(u)}\le r\quad(u\in[-1,1]^5),\qquad
 \deg P\le d,\quad P\in2^{-40}\Z[u],\quad r\in2^{-40}\Z_{\ge0}.
\]
\[
 (P,r)+(Q,s)=(P+Q,r+s),\qquad -(P,r)=(-P,r),\qquad
 \lceil a\rceil_{2^{-40}}=2^{-40}\lceil2^{40}a\rceil.
\]
For $d\in\{1,2\}$,
\[
 P=\sum_{\substack{\alpha\in\Nzero^5\\|\alpha|\le d}}p_\alpha u^\alpha,
 \qquad Q=\sum_{\substack{\beta\in\Nzero^5\\|\beta|\le d}}q_\beta u^\beta.
\]
The product enclosure is
\begin{equation}\label{gl:enclosureproduct}
 \begin{aligned}
 P'&=\sum_{\substack{\gamma\in\Nzero^5\\|\gamma|\le d}}
 \mathcal T_{2^{-40}}\!\left(
 \sum_{\substack{\alpha,\beta\in\Nzero^5\\\alpha+\beta=\gamma}}
       p_\alpha q_\beta\right)u^\gamma,\\
 r'&=\Biggl\lceil
 s\sum_{\substack{\alpha\in\Nzero^5\\|\alpha|\le d}}|p_\alpha|
 +r\sum_{\substack{\beta\in\Nzero^5\\|\beta|\le d}}|q_\beta|+rs\\
 &\qquad+
 \sum_{\substack{\alpha,\beta\in\Nzero^5\\
                  |\alpha|,|\beta|\le d,\ |\alpha|+|\beta|>d}}
         |p_\alpha q_\beta|
 \Biggr\rceil_{2^{-40}}\\
 &\qquad+2^{-40}\begin{cases}
 \#\left\{\gamma\in\Nzero^5:|\gamma|\le1,\
 \displaystyle\sum_{\substack{\alpha,\beta\in\Nzero^5\\\alpha+\beta=\gamma}}
 p_\alpha q_\beta\notin2^{-40}\Z\right\},&d=1,\\
 21,&d=2.
 \end{cases}
 \end{aligned}
\end{equation}
Since $|u^\alpha|\le1$, the discarded coefficients and input errors give
\eqref{gl:enclosureproduct}. Use rational arithmetic for each coefficient sum before truncation
and for the radius before its single ceiling.
For dyadic $a$, the degree-two rule is
\begin{equation}\label{gl:scalarmult}
 S_a(P,r)=\big(\mathcal T_{2^{-40}}(aP),
                 \lceil\abs a r\rceil_{2^{-40}}+21\cdot2^{-40}\big),
\end{equation}
$a\in\{0,1,-1\}$ acts without rounding.

For the degree-one reciprocal,
\[
 \begin{gathered}
 a=P(0),\qquad
 A=r+\sum_{\substack{\alpha\in\Nzero^5\\|\alpha|=1}}|p_\alpha|,\\
 a>A,\qquad\beta=\lceil A/a\rceil_{2^{-40}}<1.
 \end{gathered}
\]
Set $b=\mathcal T_{2^{-40}}(1/a)$. Evaluate $b-(bb)(P-a,r)$ by
\eqref{gl:enclosureproduct}, then add to its radius
\begin{equation}\label{gl:reciprocal}
 2^{-40}+\lceil2^{-40}(2b+2^{-40})A\rceil_{2^{-40}}
 +\left\lceil\frac{(b+2^{-40})\lceil\beta^2\rceil_{2^{-40}}}{1-\beta}
                                                    \right\rceil_{2^{-40}}.
\end{equation}
Indeed,
\[
 \frac1{a+x}=\frac1a-\frac{x}{a^2}+\frac{x^2}{a^2(a+x)},\qquad
 0\le a^{-1}-b<2^{-40},\qquad
 \frac{x^2}{a^2(a+x)}\le\frac{(b+2^{-40})\beta^2}{1-\beta}
 \quad(\abs x\le A).
\]

A lower endpoint is
\begin{equation}\label{gl:lower}
 \mathcal L(P,r)=p_0-r-
 \sum_{\substack{\alpha\in\Nzero^5,\ 0<|\alpha|\le d\\p_\alpha<0\text{ or }\alpha_i\text{ odd for some }1\le i\le5}}
                   \abs{p_\alpha}
\end{equation}
with
\[
 |u^\alpha|\le1,\qquad
 p_\alpha\ge0,\ \alpha\in(2\Nzero)^5
 \Longrightarrow p_\alpha u^\alpha\ge0,
 \qquad h(u)\ge\mathcal L(P,r).
\]

For polynomial trials, use degree two in \eqref{gl:linearRayleigh} and
\begin{equation}\label{gl:tau}
 \tau=2^{-18}\max\{j\in\Z:0\le j<2^{26},\ j^3<2^{54}27T^3D(c)^2\}.
\end{equation}
For $h\le40$,
\[
 c_{01}=1,\qquad c_e\in2^{-9}\Z\cap[2^{-9},1]\quad(e\ne01),
\]
so
\[
 D(c)\ge c_{01}c_{02}c_{03}\ge2^{-18},\qquad
 1<2^{54}27T^3D(c)^2,\qquad \tau>0.
\]
Use
\[
 q_e\longmapsto(q_e,0),\qquad
 q_e\longmapsto(q_e,2^{-39})\quad(\mathrm{FCC}),\qquad
 a\in\mathbb Q\longmapsto(\mathcal T_{2^{-40}}(a),2^{-40}).
\]
Compute in order
\[
 \begin{gathered}
 \sum_{0\le i<j\le3}S_{C_{ij}(c)}q_{ij},\qquad
 \prod_{e\in J}q_e=(q_iq_j)q_k\quad(J=\{i,j,k\},\ i<j<k),\\
 S_\tau M-\left(\sum_{0\le i<j\le3}S_{C_{ij}(c)}q_{ij}\right)K.
 \end{gathered}
\]
Products use \eqref{gl:enclosureproduct}; dyadic scalars use $S_a$.

At each cubature node,
\[
 \begin{aligned}
 x_{J,i}&\longmapsto(\mathcal T_{2^{-40}} x_{J,i},6\cdot2^{-40})
 &&(1\le i\le3),\\
 \Gamma_{ij}&=\sum_{0\le a<b\le3}S_{(b_a-b_b)_i(b_a-b_b)_j/16}q_{ab}
 &&(1\le i,j\le3).
 \end{aligned}
\]
For FCC, apply in order
\[
 (z_1,z_2,z_3)\leftarrow(x_3,x_1+x_2,x_1-x_2),
\]
\[
 \begin{aligned}
 (\Gamma_{1j},\Gamma_{2j},\Gamma_{3j})
 &\leftarrow(\Gamma_{3j},\Gamma_{1j}+\Gamma_{2j},\Gamma_{1j}-\Gamma_{2j})
 &&(1\le j\le3),\\
 (\Gamma_{i1},\Gamma_{i2},\Gamma_{i3})
 &\leftarrow(\Gamma_{i3},\Gamma_{i1}+\Gamma_{i2},\Gamma_{i1}-\Gamma_{i2})
 &&(1\le i\le3).
 \end{aligned}
\]
Each assignment is simultaneous. Apply \eqref{gl:scalarmult} with factor
\[
 2^{-24}\left\lfloor
 \frac{12\cdot2^{24}}{3(c_{01}+c_{02}+c_{13}+c_{23})+4(c_{03}+c_{12})}
 \right\rfloor
\]
once to each $z_i$ and twice successively to each $\Gamma_{ij}$;
evaluate the scale once and reuse it.
For cubic monomials and their derivatives, start with the enclosure $(1,0)$ and use successive powers,
grouped as $(x_1^ax_2^b)x_3^c$. For quintics, in increasing degree,
\[
 x^0=(1,0),\qquad x^{e_i}=x_i,\qquad
 x^\alpha=x^{\alpha-e_i}x_i,\quad
 i=\min\{j\in\{1,2,3\}:\alpha_j>0\},\quad |\alpha|\ge2,
\]
For FCC, use coordinate powers grouped as
\[
 \text{value},\ \partial_1:\ z_1^a(z_2^bz_3^c),\qquad
 \partial_2,\ \partial_3:\ (z_1^az_2^b)z_3^c,
\]
Lower the differentiated exponent; omit zero coefficients.
For independently evaluated enclosures of $f$ and $\partial_i f$, each tree gives
\[
 \begin{split}
 &\text{mass: }\quad\sum_{\substack{y\in[-1/2,1/2]^3\\w_y\ne0}} w_y(ff),\\
 &\text{energy: }\quad\sum_{\substack{y\in[-1/2,1/2]^3\\w_y\ne0}} w_y
       \sum_{i=1}^3\sum_{j=1}^{i}
                S_{2-\mathbf1_{i=j}}\big(\Gamma_{ij}((\partial_i f)(\partial_j f))\big).
 \end{split}
\]
Omit zero factors; multiply by $(q_iq_j)q_k$, then sum over trees.

For $A=(0,a_1,a_2,a_3)\in\mathcal T_+$ and $s\in\{-1,1\}^3$, let
$V_q(w)$ denote the corresponding vertex or face average in \eqref{gl:vertex}:
\[
 H_{ij}=\sum_{k=1}^3
 s_i\bigl(V_q(\diag(s_1,s_2,s_3)a_k)\bigr)_i
 \bigl([a_1,a_2,a_3]^{-1}\bigr)_{kj}\qquad(1\le i,j\le3).
\]
Use affine coefficient multiplication and addition with all radii, and
\[
 a\in2^{-40}\Z:\ (a,0),\qquad
 a\notin2^{-40}\Z:\ (\mathcal T_{2^{-40}}(a),2^{-40}).
\]
With indices modulo three,
\[
 \begin{aligned}
 (\adj H)_{ij}&=H_{j+1,i+1}H_{j+2,i+2}-H_{j+1,i+2}H_{j+2,i+1},\\
 \det H&=\sum_{i=1}^3H_{1i}(\adj H)_{i1}.
 \end{aligned}
\]
Form $(s_i s_j\Gamma_{ij})_{i,j=1}^3$ by sign changes without rounding.
Multiply first on the left by $\adj H$, then on the right by $(\adj H)^T$;
reuse the resulting entries below. Each of the $378$ reference integrals below, evaluated by
\eqref{lo:moments}, is enclosed once by its nearest $2^{-40}$-multiple
(ties even), with radius $2^{-40}$; the formulas below reuse these enclosures:
\begin{equation}\label{gl:transportenclosure}
 \begin{aligned}
 M&=\sum_{1\le i\le j\le3}(2-\mathbf1_{i=j})v_iv_j
   \sum_{A\in\mathcal T_+}\sum_{s\in\{-1,1\}^3}s_is_j
       \left((\det H)\int_A P_iP_j\right),\\
 K&=\sum_{1\le i\le j\le3}(2-\mathbf1_{i=j})v_iv_j
   \sum_{A\in\mathcal T_+}\sum_{s\in\{-1,1\}^3}s_is_j\\
 &\qquad\times\Biggl(\Biggl(
       \sum_{1\le k\le l\le3}
       \Bigl((\adj H)(s_a s_b\Gamma_{ab})_{a,b=1}^3(\adj H)^T\Bigr)_{kl}\\
 &\hspace{100pt}{}\times\int_A(\partial_kP_i\partial_lP_j
       +\mathbf1_{k<l}\partial_lP_i\partial_kP_j)\Biggr)(\det H)^{-1}\Biggr).
 \end{aligned}
\end{equation}
Keep $s_is_j\in\{-1,1\}$ unrounded and enclose $(2-\mathbf1_{i=j})v_iv_j$ as one rational.
Use the affine reciprocal for $\det H$ and
\[
 D=\sum_{\substack{\{a,b,c\}\in\{J_1,\ldots,J_{16}\}\\a<b<c}}(q_aq_b)q_c.
\]
Enclose $T$ as one rational constant; evaluate $TM$ once and reuse it.
Test the lower endpoints of $M$ and
\begin{equation}\label{gl:cubicexpression}
 ((TM)(TM))(TM)-((DK)K)K.
\end{equation}

For center selection, round rational inputs once and every operation and
square root to binary64 (ties even, signed zero). Sums start at $0$ in
displayed order, otherwise increasing-index order; square roots use squared dyadic midpoints.

For cubic and quintic trials, $s=3,5$, let $m_\alpha$ be the rounded
cubature moment of $x^\alpha$ at $c$, $\alpha\in\Nzero^3$, $|\alpha|\le2s$
even. Use lexicographic tree and node-coordinate order, rounded coordinates
$x_J(y;c)$, successive powers, and weights
\[
 \left[\left(\prod_{e\in J}c_e\right)\mathcal T_{2^{-40}}(w_y)\right],
\]
with rational $w_y$. Accumulate $\alpha=(a,b,d)$ terms as
$((wx_1^a)x_2^b)x_3^d$. The lower triangles are
\begin{equation}\label{gl:gram}
 M_{\alpha\beta}=m_{\alpha+\beta},\qquad
 K_{\alpha\beta}=\sum_{\substack{1\le i,j\le3\\\alpha_i\beta_j>0}}
       (\alpha_i\beta_j\Gamma_{ij}(c))m_{\alpha+\beta-e_i-e_j}.
\end{equation}
$i$ is the outer summation index; complete by symmetry.

For FCC at $q=c$,
\[
 (z_1,z_2,z_3)\leftarrow([x_3],[x_1+x_2],[x_1-x_2]).
\]
Use the displayed row/column updates and one coordinate or two successive
metric multiplications by the factor in \eqref{gl:eta}. Evaluate it once;
round every multiplication to binary64.
Use $g_i$ without coefficient truncation, in the monomial order
\[
 z_1,z_1z_2^2,z_1z_3^2,z_1z_2^4,z_1z_3^4,z_1z_2^2z_3^2,
 z_1^3,z_1^3z_2^2,z_1^3z_3^2,z_1^5.
\]
Preserve this order under coordinate interchange. In the lower triangle
of \eqref{gl:gram}, multiply by both monomial coefficients and sum pairs
lexicographically. For transport, use the constant coefficients of the six
affine enclosures of each of $M(P_i,P_j)$ and $K(P_i,P_j)$.

For each $(M,K)$, the serial Cholesky recurrence is
\[
 \begin{split}
 r_0&=M_{ij},\qquad
 r_k=[r_{k-1}-[L_{ik}L_{jk}]]\quad(1\le k<j\le i),\\
 L_{ij}&=
 \begin{cases}
 [\sqrt{r_{j-1}}],&i=j,\\
 [r_{j-1}/L_{jj}],&i>j.
 \end{cases}
 \end{split}
\]
A nonpositive diagonal pivot rejects the test. By forward substitution,
$W=L^{-1}$; compute $WK$, then the lower triangle of $E=(WK)W^T$, and
complete by symmetry.

Start $V=I$. For twelve lexicographic sweeps of $i<j$, update whenever
$|E_{ij}|\ge[10^{-13}]$:
\begin{equation}\label{gl:jacobi}
 \begin{gathered}
 u=\frac{E_{jj}-E_{ii}}{2E_{ij}},\qquad
 t=\begin{cases}
 \dfrac1{\abs u+\sqrt{1+u^2}},&u>0\text{ or }u\text{ is }+0,\\[4pt]
 \dfrac{-1}{\abs u+\sqrt{1+u^2}},&u<0\text{ or }u\text{ is }-0,
 \end{cases}\\
 a=\frac1{\sqrt{1+t^2}},\qquad b=at.
 \end{gathered}
\end{equation}
Here $+0$ and $-0$ are the signed binary64 zeros. Apply
\[
 \begin{split}
 &(E_{ii},E_{jj})\leftarrow(E_{ii}-tE_{ij},E_{jj}+tE_{ij}),\qquad
 E_{ij}=E_{ji}\leftarrow0,\\
 &(E_{ki},E_{kj})=(x,y)\leftarrow(ax-by,bx+ay)
       \quad(k\ne i,j),\\
 &E_{ik}\leftarrow E_{ki},\qquad E_{jk}\leftarrow E_{kj},
 \end{split}
\]
and the same two-column update rowwise to $V$. For the first minimizing
diagonal index $i$,
\[
 V_{\cdot i}\leftarrow W^TV_{\cdot i},\qquad
 v_j=2^{-24}\operatorname{nint}\!\left(
 2^{24}\left[\frac{V_{ji}}{\max_k|V_{ki}|}\right]\right).
\]
Fix the maximum before normalization; reject the test if $v=0$.

For $|w|=h$, stop at the first successful test in this order:
\begin{enumerate}
 \item Chamber $\longrightarrow$ corner $\longrightarrow$ local
 (Section~\ref{sec:global}).
 \item
 \[
 h\ge12,\quad q_e^->0\ (e\in E),\quad
 \min_{0\le i<j\le3}c_{ij}\ge\tfrac12
 \quad\Longrightarrow\quad\text{transported test}.
 \]
 \[
 h\ge11,\quad w_1=\cdots=w_{11}=1
 \quad\Longrightarrow\quad\text{omit subsequent tests}.
 \]
 \item
 \[
 h\ge12,\quad(h\text{ even or }h\ge35)
 \quad\Longrightarrow\quad\text{cubic test}.
 \]
 \item For $h\ge20$ with $h$ even or $h\ge35$,
 \[
 \min(c_{03},c_{12})\le\tfrac13,\quad
 \min(c_{02},c_{13},c_{23})\ge\tfrac23
 \quad\Longrightarrow\quad\text{FCC test},
 \]
 followed by the quintic test.
\end{enumerate}
A failed reciprocal condition rejects that trial. Use the preceding center
trials and \eqref{gl:recursion} through depth $41$;
\[
 |w|\le40\quad\Longrightarrow\quad q_e^\pm\in2^{-8}\Z.
\]

The first-success counts are
\[
 \begin{aligned}
 (\text{chamber},\text{corner},\text{local})&=(8655,1813,507),\\
 (\text{cubic},\text{quintic},\text{FCC},\text{transport})
   &=(14556,3373,1653,61342).
 \end{aligned}
\]
The minimum lower endpoints are
\[
 \begin{aligned}
 \min\mathcal L(\tau M-\left(\sum_{0\le i<j\le3}S_{C_{ij}(c)}q_{ij}\right)K)_{\rm cubic}&=12339663\cdot2^{-40}>0,\\
 \min\mathcal L(\tau M-\left(\sum_{0\le i<j\le3}S_{C_{ij}(c)}q_{ij}\right)K)_{\rm quintic}&=2997515\cdot2^{-40}>0,\\
 \min\mathcal L(\tau M-\left(\sum_{0\le i<j\le3}S_{C_{ij}(c)}q_{ij}\right)K)_{\rm FCC}&=3874396\cdot2^{-40}>0,\\
 \min\mathcal L(T^3M^3-DK^3)_{\rm transport}&=799999930\cdot2^{-40}>0,\\
 \min\mathcal L(M)_{\rm transport}&=684637240747\cdot2^{-40}>0.
 \end{aligned}
\]
Finally,
\[
 \#\mathcal W_0=1,\qquad
 \#\mathcal W_{h+1}
 =2\#\{w\in\mathcal W_h:B_w\text{ is not accepted}\},\qquad
 \#\mathcal W_{37}=4.
\]
\[
 \begin{aligned}
 \sum_{h=0}^{37}\#\{w\in\mathcal W_h:B_w\text{ is accepted}\}&=91899,\\
 \sum_{h=0}^{36}\#\{w\in\mathcal W_h:B_w\text{ is not accepted}\}&=91898,
 \end{aligned}
\]
\[
 \{w\in\mathcal W_{37}:B_w\text{ is not accepted}\}=\varnothing
 \quad\Longrightarrow\quad\mathcal W_{38}=\varnothing.
\]
\end{proof}

\bigskip
\noindent
{\bf Acknowledgements.}
The authors used AI to assist the long but finite computations.

\end{document}